\documentclass[12pt]{amsart} \usepackage{amsmath,centernot,amssymb,leftindex}
\usepackage{enumitem}
\setlist[enumerate]{font={\textrm}}
\setlist[enumerate,1]{label={(\roman*)}}
\usepackage[backend=biber,style=alphabetic,sorting=nty,sortcites]{biblatex}
\DeclareFieldFormat{url}{%
	\iffieldundef{doi}{%
		\mkbibacro{URL}\addcolon\space\url{#1}%
	}{%
	}%
}

\DeclareFieldFormat{urldate}{%
	\iffieldundef{doi}{%
		\mkbibparens{\bibstring{urlseen}\space#1}%
	}{%
	}%
}
\usepackage[hidelinks]{hyperref}

\newcommand{\Bohr}{\operatorname{Bohr}}

\newcommand{\mb}{\mathbf}

\newtheorem{theorem}{Theorem}
\newtheorem{lemma}[theorem]{Lemma}

\newtheorem{corollary}[theorem]{Corollary}
\newtheorem*{claim}{Claim}
\numberwithin{theorem}{section}
\numberwithin{equation}{section}

\theoremstyle{definition}

\newtheorem{observation}[theorem]{Observation}
\newtheorem{definition}[theorem]{Definition}
\newtheorem{question}[theorem]{Question}

\newtheorem{remark}[theorem]{Remark}
\newtheorem{example}[theorem]{Example}

\title[Sumsets with one large summand]{Discrete sumsets with one large summand}

\date{\today}

\author{John T. Griesmer}
\email{jtgriesmer@gmail.com}
\address{Department of Applied Mathematics and Statistics, Colorado School of Mines, Golden, Colorado}
\usepackage{amstext}

\begin{document}

	\begin{abstract}
		If $A$ and $B$ are subsets of an abelian group, their sumset is $A+B:=\{a+b:a\in A, b\in B\}$.  We study sumsets in discrete abelian groups, where at least one summand has positive upper Banach density. 
		
		Jin proved in \cite{Jin_SumsetPhenom} that if $A$ and $B$ are sets of integers having positive upper Banach density, then $A+B$ is piecewise syndetic.  Bergelson, Furstenberg, and Weiss \cite{BergelsonFurstenbergWeiss} improved the conclusion to ``$A+B$ is piecewise Bohr.''  In \cite{BeiglbockBergelsonFish} this was shown to be qualitatively optimal, in the sense that if $C\subseteq \mathbb Z$ is piecewise Bohr, then there are $A, B\subseteq \mathbb Z$ having positive upper Banach density such that $A+B\subseteq C$.
		
		We improve these results by establishing a strong correspondence between sumsets in discrete abelian groups, level sets of convolutions in compact abelian groups, and sumsets in compact abelian groups. Our proofs avoid measure preserving dynamics and nonstandard analysis, and our results apply to discrete abelian groups of any cardinality.
	\end{abstract}
	\maketitle
	
	\setcounter{tocdepth}{1}
	\tableofcontents
	
	\section{The Steinhaus lemma and its density analogues}
	
	Let $G$ be an abelian group. If $A$ and $B$ are subsets of $G$, their \emph{sumset} is $A+B:=\{a+b:a\in A, b\in B\}$.  If $t\in \Gamma$, write $A+t$ for the \emph{translate} $\{a+t:a\in A\}$.  We are interested in the structure of $A+B$ under various hypotheses saying ``$A$ and  $B$ are large.''  
	
	\subsection{The Steinhaus lemma}\label{sec:SteinhausReview}
	
	H.~Steinhaus \cite[Th{\'e}or\`eme VII]{Steinhaus_SurLesDistances} proved that if $A, B\subseteq \mathbb R$ both have positive Lebesgue measure, then $A+B$ contains an interval.  The relevant properties of Lebesgue measure are  retained by Haar measure on locally compact abelian (LCA) groups, leading to the following generalization.  
	\begin{theorem}\label{th:SteinhausGeneral}
		Let $K$ be an LCA group with Haar measure $\mu$.  If $\mu(A)>0$ and $\mu(B)>0$, then $A+B$ contains a nonempty open subset of $K$.
	\end{theorem}
We will prove a special case of Theorem \ref{th:SteinhausGeneral} in \S\ref{sec:SteinhausProof}. Weil \cite[p.50]{Weil2ndEd} proves the above for arbitrary locally compact (including nonabelian) groups.   See \cite{Stromberg} for a short proof of the special case where $B=A^{-1}$.
	
	Our main results are Theorem  \ref{th:Main} and Theorem \ref{th:MainGeneral}; these are  analogues of  Theorem \ref{th:SteinhausGeneral} where Haar measure is replaced by a finitely additive translation-invariant measure.  Such measures are usually discussed in terms of invariant means on the space of bounded functions on the underlying group.
		\subsection{Means on \texorpdfstring{$\ell^{\infty}(\Gamma)$}{l infty Gamma}}  Let $\Gamma$ be a discrete abelian group. We write $\ell^{\infty}(\Gamma)$ for the Banach space of bounded complex-valued functions on $\Gamma$, equipped with the supremum norm $\|\cdot\|_{\infty}$.  For a function $f$ defined on $\Gamma$ and $\gamma \in \Gamma$, we write $f_\gamma$ for the \emph{translate} given by $f_\gamma(x):=f(x-\gamma)$.  A \emph{mean} on $\ell^\infty(\Gamma)$ is a linear functional $m:\ell^\infty(\Gamma)\to \mathbb C$ satisfying $m(1_\Gamma) = 1$ and
	$m(f)\geq 0$ if $f(x)\geq 0$ for all $x\in \Gamma$. 	
	
		We write $f\sim_{m} g$ if $m(|f-g|)=0$. 
	
	Note that each mean $m$ induces a \emph{finitely additive} probability measure $m^*$ on the algebra of subsets of $\Gamma$, given by $m^*(A)=m(1_A)$.  
	
	We write $A\subset_m B$ if $m^*(A\setminus B)=0$, and we write $A\sim_{m} B$ if $m^{*}(A\triangle B)=0$, or equivalently, if $1_{A} \sim_{m} 1_{B}$.  
	
	We will abuse notation and write $m(A)$ for $m^{*}(A)$.
	
	If a mean $m$ satisfies
	$m(f_\gamma)=m(f)$ for all $f\in \ell^{\infty}(\Gamma)$ and $\gamma \in \Gamma$,	we say that $m$ is an \emph{invariant mean}; invariant means are sometimes called ``Banach means'' or ``Banach mean values'' in the literature, cf.~\cite{Folner_Generalization}, \cite{Folner_NoteOnGeneralization}.

The set $\mathcal M(\Gamma)$ of means on $\ell^{\infty}(\Gamma)$ is a weak$^{*}$-compact convex subset of $\ell^{\infty}(\Gamma)^{*}$. Translation on $\Gamma$ induces an action $\tau$ of $\Gamma$ on $\mathcal M(\Gamma)$ by continuous linear operators: $(\tau_{\gamma} m) (f):=m(f_{\gamma})$.  The Markov-Kakutani fixed point theorem \cite[p.~151]{Conway_FunctionalAnalysis} guarantees that there is at least one mean $m\in \mathcal M(\Gamma)$ satisfying $\tau_{\gamma}m=m$ for all $\gamma\in \Gamma$.  In other words, there is at least one invariant mean on $\ell^{\infty}(\Gamma)$. 

We let $\mathcal M_{\tau}(\Gamma)$ denote the set of invariant means on $\ell^\infty(\Gamma)$.   It is straightforward to verify that $\mathcal M_{\tau}(\Gamma)$ is a weak$^*$-compact convex subset of $\ell^\infty(\Gamma)^*$. The Krein-Milman theorem \cite[p.~142]{Conway_FunctionalAnalysis} then implies that the set $\mathcal M_{\tau}^{ext}(\Gamma)$ of extreme points of $\mathcal M_{\tau}(\Gamma)$ is nonempty.  We call  elements of $\mathcal M_{\tau}^{ext}(\Gamma)$ \emph{extreme} invariant means.

	\begin{definition}\label{def:UBD}
		Let $A\subseteq \Gamma$.  The \emph{upper Banach density} of $A$ is 
		\[d^*(A):=\sup\{m(A):m \in \mathcal M_{\tau}(\Gamma)\}.\]
	\end{definition}
See Lemma \ref{lem:UBDequivalents} for conditions equivalent to $A$ having positive upper Banach density.
	\begin{observation}\label{obs:UBDrealized} For every $A\subseteq \Gamma$, there is an extreme invariant mean $\nu$ such that $d^{*}(A)=\nu(1_{A})$. To see this, note that $\lambda \mapsto \lambda(1_{A})$ is a continuous linear functional on $\ell^{\infty}(\Gamma)^{*}$ and $\mathcal M_{\tau}(\Gamma)$ is a compact convex subset thereof. Proposition 7.9 in \cite[p.~144]{Conway_FunctionalAnalysis} then implies that the supremum $\{m(A):m\in \mathcal M_{\tau}(\Gamma)\}$ is attained by at least one extreme point of $\mathcal M_{\tau}(\Gamma)$. See Lemma \ref{lem:UBDrealized} below for details. 
	\end{observation}

\subsection{Characters} 
	
	We use standard definitions and background presented in e.g.~Chapter 4 of \cite{Folland_CourseInAHA} or Chapters 1 and 2 of \cite{Rudin_Fourier}.  
	
	Write $\mathcal S^1$ for the group $\{z\in \mathbb C: |z|=1\}$ with the group operation of multiplication and the usual topology. 
	
	Let $G$ be an LCA group. A \emph{character} of $G$ is a continuous homomorphism $\chi:G\to \mathcal S^{1}$; the set of all characters of $G$ is denoted $\widehat{G}$.  The \emph{trivial character} $\chi_{0}$ is constant: $\chi_{0}(\gamma)=1$ for all $\gamma\in \Gamma$. With the group operation of pointwise multiplication and the topology of uniform convergence on compact subsets, $\widehat{G}$ is an LCA group.  When $G$ is compact, $\widehat{G}$ is discrete, and when $G$ is discrete, $\widehat{G}$ is compact.  A \emph{trigonometric polynomial} is a finite linear combination of characters of $G$.  A function $f:G\to \mathbb C$ is \emph{uniformly almost periodic} if it is a uniform limit of trigonometric polynomials.

	\subsection{Bohr compactification}\label{sec:BohrCompactification}
Let $\Gamma$ be a discrete abelian group. The \emph{Bohr compactification of $\Gamma$} is a compact Hausdorff abelian group  $b\Gamma$ together with a one-to-one homomorphism $\iota:\Gamma \to b\Gamma$ such that
\begin{enumerate}
	\item[(i)] $\iota(\Gamma)$ is topologically dense in $b\Gamma$;
	
	\item[(ii)] for all $\chi\in \widehat{\Gamma}$, there is a character $\tilde{\chi}\in \widehat{b\Gamma}$ such that $\tilde{\chi}\circ \iota = \chi$.
\end{enumerate}
We will identify $\Gamma$ with its image $\iota(\Gamma)$, and consider $\Gamma$ as a topologically dense subgroup of $b\Gamma$, so that each $\chi\in \widehat{\Gamma}$ is the restriction of a continuous $\tilde{\chi}\in \widehat{b\Gamma}$: $\chi = \tilde{\chi}|_\Gamma$.  

The group $b\Gamma$ may be constructed as the Pontryagin dual of $\widehat{\Gamma}_d$, where $\widehat{\Gamma}_d$ is the group $\widehat{\Gamma}$ with the discrete topology instead of the usual topology.  In this construction, the embedding map $\iota$ is given by $\iota(\gamma):=e_\gamma$, where $e_\gamma(\chi)=\chi(\gamma)$; see Section 1.8 of \cite{Rudin_Fourier} or Section 4.7 of \cite{Folland_CourseInAHA} for details.

The following is the special case of \cite[ Theorem 4.79]{Folland_CourseInAHA} for discrete abelian groups.

\begin{theorem}\label{th:BohrCompactificationAndAP} 
	If $f$ is a bounded  function on $\Gamma$, the following are equivalent: 
	\begin{enumerate}
		\item[a.]  $f$ is the restriction to $\Gamma$ of a continuous function on $b\Gamma$.
		\item[b.]  $f$ is a uniform limit of linear combinations of characters on $\Gamma$.
		\item[c.]  $f$ is uniformly almost periodic.	
	\end{enumerate}
\end{theorem}
We write $\mu_{b\Gamma}$ for Haar probability measure on $b\Gamma$.

\subsection{Main Theorem}	
	The next theorem is our main result; it relates sumsets in   discrete abelian groups to sumsets in compact abelian groups.  In \S\ref{sec:SumsetPhenomena} we show how it extends the results mentioned in the abstract.  Here ``$F_{\sigma}$'' means ``countable union of compact sets.''  Note that when $A$ and $B$ are $F_{\sigma}$ subsets of $b\Gamma$, $A+B$ is also $F_{\sigma}$, implying that $A+B$ is $\mu_{b\Gamma}$-measurable.
	
	\begin{theorem}\label{th:Main}
	Let $\Gamma$ be a discrete abelian group.  Let $\nu$ be an extreme invariant mean on $\ell^{\infty}(\Gamma)$, $m$ an invariant mean on $\ell^{\infty}(\Gamma)$, and $A$, $B\subseteq \Gamma$.  
	There are $F_{\sigma}$ sets $\tilde{A}_{\nu}$, $\tilde{B}_{m}\subseteq b\Gamma$  such that $\mu_{b\Gamma}(\tilde{A}_{\nu})\geq \nu(A)$, $\mu_{b\Gamma}(\tilde{B}_{m})\geq m(B)$, and
	\begin{equation}
		\nu(A+B)\geq \mu_{b\Gamma}(\tilde{A}_{\nu}+\tilde{B}_{m}).
	\end{equation}
	Furthermore,
		\begin{enumerate}
		\item\label{item:AVBV} if $V\subseteq b\Gamma$ is compact, then 
		\[\mu_{b\Gamma}(\tilde{A}_{\nu}\cap V)\geq \nu(A\cap V) \quad \text{and} \quad \mu_{b\Gamma}(\tilde{B}_{m}\cap V)\geq m(B\cap V);\]
		
						\item\label{item:VbGammaClopen} if $V\subseteq b\Gamma$ is clopen, then
	\[\mu_{b\Gamma}(\tilde{A}_{\nu}\cap V)=0\iff \nu(A\cap V)=0 \quad \text{and} \quad \mu_{b\Gamma}(\tilde{B}_{m}\cap V)=0\iff m(B\cap V)=0;\]
		
		\item\label{item:VcapGamma}  if $V\subseteq \tilde{A}_{\nu}+\tilde{B}_{m}$ is compact, then $V\cap \Gamma \subset_{\nu} A+B$.


	\end{enumerate}
	\end{theorem}
	We will prove Theorem \ref{th:Main} as a special case of Theorem \ref{th:MainGeneral} below; the latter obtains the same conclusion under a weaker hypothesis on $B$.

		\subsection{Large sets}\label{sec:JinBFW}   We list some common definitions for large  subsets of discrete abelian groups.  Equivalent formulations are provided in \S\ref{sec:LargeSubsets}.
	
	Let $\Gamma$ be an arbitrary discrete abelian group (countable or uncountable) and $\gamma_{0}\in \Gamma$.   We say $U$ is a \emph{Bohr neighborhood of $\gamma_{0}$} if there is an open subset $V\subseteq b\Gamma$ such that $\gamma_{0}\in V$ and $V\cap \Gamma \subseteq U$.  By Lemma \ref{lem:BohrEquivalents}, $U$ is a Bohr neighborhood of $\gamma_{0}$ if and only if there is a uniformly almost periodic function $\phi$ supported on $U$ such that
	 $\phi(\gamma_{0})\neq 0$.

	A set $T\subseteq \Gamma$ is \emph{thick} if for all finite $F\subseteq \Gamma$, there is a $t\in \Gamma$ such that $F+t\subseteq T$.   By Lemma \ref{lem:ThickEquivalents}, $T$ is thick if and only if there is an invariant mean $m$ on $\ell^{\infty}(\Gamma)$ such that $m(T)=1$.
	
	A set $S\subseteq \Gamma$ is \emph{syndetic} if there is a finite set $F\subseteq \Gamma$ such that $S+F=\Gamma$. By Lemma \ref{lem:syndeticEquivalents}, $S$ is syndetic if and only if $S$ has nonempty intersection with every thick set.  Every Bohr neighborhood is syndetic -- see \S\ref{sec:BohrNhoods}.
	
	We say that $C\subseteq \Gamma$ is \emph{piecewise syndetic} if there is an $m\in \mathcal M_{\tau}(\Gamma)$ and a syndetic set $S$ such that $S\subset_{m} C$.  Equivalently (Lemma \ref{lem:PWSequivalents}),  $C$ is piecewise syndetic if $C$ contains a set of the form $S\cap T$, where $S$ is syndetic and $T$ is thick. 
	
	We say that $C$ is \emph{piecewise Bohr} if there is an $m\in \mathcal M_{\tau}(\Gamma)$ and a Bohr neighborhood $U\subseteq \Gamma$ such that $U\subset_{m}C$.  Equivalently (Lemma \ref{lem:PWBohrEquivalents}), $C$ is piecewise Bohr if there is a thick set $T$ and a Bohr neighborhood $U$ of some $\gamma\in \Gamma$ such that $T\cap U\subseteq \Gamma$.
	
	Every piecewise Bohr set is piecewise syndetic, but some syndetic sets are not piecewise Bohr;
	examples are given in \cite[Theorem 4.3]{BergelsonFurstenbergWeiss} and in Section 2.5 of \cite{Griesmer_SumsetsDenseSparse}.

	\subsection{Sumset phenomena}\label{sec:SumsetPhenomena}	  Jin \cite{Jin_SumsetPhenom} used nonstandard analysis to give a new proof of the Steinhaus lemma in $\mathbb R$, and to prove a new result about sumsets in $\mathbb Z$: Corollary 3 of \cite{Jin_SumsetPhenom} says that if $A, B\subseteq \mathbb Z$ both have positive upper Banach density, then $A+B$ is piecewise syndetic.  Under the same hypothesis, Bergelson, Furstenberg, and Weiss  improved the conclusion to say in  that $A+B$ is piecewise Bohr; this is  \cite[Theorem I]{BergelsonFurstenbergWeiss}.  The latter result was extended to countable amenable groups by Beiglb\"ock, Bergelson, and Fish in \cite{BeiglbockBergelsonFish}, and to arbitrary countable groups by Bj\"orklund and Fish in \cite{BjorklundFish_ProductSet} and by Bj\"orklund in \cite{Bjorklund_ProductSetMeasured}.   As Corollary \ref{cor:PWB} shows, Theorem \ref{th:Main}  recovers the abelian cases of these results and extends them to discrete abelian groups of arbitrary cardinality.   
	
	In \cite[Theorem 1.4]{Griesmer_SumsetsDenseSparse} we showed that when $d^{*}(A)>0$ and $B$ satisfies a weaker positive density condition (cf.~\S\ref{sec:FSmeans} below), we may still conclude that $A+B$ is piecewise Bohr.  Corollary \ref{cor:Recover} extends this result to arbitrary discrete abelian groups. 
	 
	\begin{corollary}\label{cor:PWB}
		If $\Gamma$ is a discrete abelian group and $A, B\subseteq \Gamma$ have positive upper Banach density, then $A+B$ is piecewise Bohr (and therefore piecewise syndetic).  Furthermore, if $\nu\in \mathcal M_{\tau}^{ext}(\Gamma)$ and $\nu(A)>0$, then there is a Bohr neighborhood $U$ such that $U\subset_{\nu} A+B$.
	\end{corollary}

\begin{proof}
	Assuming $A$, $B\subseteq \Gamma$ have positive upper Banach density, apply Observation \ref{obs:UBDrealized} to find $\nu$, $m\in \mathcal M_{\tau}^{ext}(\Gamma)$ such that $\nu(A)=d^{*}(A)>0$ and $m(B)=d^{*}(B)>0$.  Let $\tilde{A}_{\nu}, \tilde{B}_{m}\subseteq b\Gamma$ be as in Theorem \ref{th:Main}, so that $\mu_{b\Gamma}(\tilde{A})$, $\mu_{b\Gamma}(\tilde{B})>0$.  By Theorem \ref{th:SteinhausGeneral}, $\tilde{A}_{\nu}+\tilde{B}_{m}$ contains a nonempty open set $\tilde{V}$.  Since $b\Gamma$ is compact Hausdorff, there is a compact subset $\tilde{V}'\subseteq \tilde{V}$ with nonempty interior.  Since $\tilde{V}'$ is a compact subset of $\tilde{A}_{\nu}+\tilde{B}_{m}$, Theorem \ref{th:Main} says that $\tilde{V}'\cap \Gamma \subset_{\nu} A+B$.   Since $\tilde{V}'$ has nonempty interior, $U:=\tilde{V}'\cap \Gamma$ is a Bohr neighborhood satisfying $U\subset_{\nu}A+B$.  Thus $A+B$ is piecewise Bohr.
\end{proof}

\begin{remark}	In Corollary \ref{cor:PWB} the conclusion that $A+B$ is piecewise Bohr cannot be improved: Theorem 4 of \cite{BeiglbockBergelsonFish} says that if $\Gamma$ is a countable abelian group and $C\subseteq \Gamma$ is piecewise Bohr, then there are sets $A, B\subseteq \Gamma$ such that $d^*(A)>0$, $d^*(B)>0$, and $A+B\subseteq C$. 
\end{remark}

\begin{remark}
	Example \ref{ex:NeedExtreme} below shows  that the hypothesis $\nu\in \mathcal M_{\tau}^{ext}(\Gamma)$ in  Corollary \ref{cor:PWB} cannot be weakened to $\nu \in \mathcal M_{\tau}(\Gamma)$.
\end{remark}

	\subsection{F{\o}lner sequences and F{\o}lner nets}
	
	Let $\mathbf F:=(F_n)_{n\in \mathbb N}$ be a sequence of finite subsets of $\Gamma$.  For $A\subseteq \Gamma$, we write
	\[
	\underline{d}_{\mb F}(A):= \liminf_{n\to\infty}\frac{|A\cap F_n|}{|F_n|}, \qquad 	\bar{d}_{\mb F}(A):= \limsup_{n\to\infty}\frac{|A\cap F_n|}{|F_n|}
	\]
	for the \emph{lower density of $A$ with respect to $\mb F$} and \emph{upper density of $A$ with respect to $\mb F$}, respectively.  
	 If $\underline{d}_{\mb F}(A)=\bar{d}_{\mb F}(A)$ we write $d_{\mb F}(A)$ for the common value, which we call the \emph{density of $A$ with respect to $\mathbf F$}.
	
	If $\mathbf F:=(F_j)_{j\in I}$ is a net of finite subsets of $\Gamma$ indexed by a directed set $I$, we define the lower and upper density of $A$ with respect to $\mathbf F$ to be $\underline{d}_{\mb F}(A):=\liminf_{j\in I} \frac{|A\cap F_j|}{|F_j|}$, and $\bar{d}_{\mb F}(A):=\limsup_{j\in I} \frac{|A\cap F_j|}{|F_j|}$, respectively.  When $f\in \ell^{\infty}(\mathbb Z)$, we write
	\[
	\underline{\mathbb E}_{\mb F}(f):=\liminf_{j\in I} \frac{1}{|F_{j}|}\sum_{\gamma\in F_{j}}f(\gamma)
	\]
	and $\bar{\mathbb E}_{\mb F}(f)$ for the corresponding $\limsup$.  When $\underline{\mathbb E}_{\mb F}(f)=\bar{\mathbb E}_{\mb F}(f)$, we write $\mathbb E_{\mb F}(f)$  for the common value.
	
	We say that $\mb F$ is a \emph{F{\o}lner sequence} if $\lim_{n\to\infty} \frac{|F_n\triangle (F_n+\gamma)|}{|F_n|}=0$ for all $\gamma\in \Gamma$.  Likewise, we say $(F_j)_{j\in I}$ is a \emph{F{\o}lner net} if $\lim_{j\in I} \frac{|F_j\triangle (F_j+\gamma)|}{|F_j|}=0$ for all $\gamma\in \Gamma$.  
	
	It is straightforward to verify that when $\mb F$ is F{\o}lner sequence (or F{\o}lner net), we have 
\begin{equation}\label{eq:FolnerTransInv}
	\underline{d}_{\mb F}(A+t)=\underline{d}_{\mb F}(A) \quad \text{ and }\quad  \bar{d}_{\mb F}(A+t)=\bar{d}_{\mb F}(A) \text{ for all } A\subseteq \Gamma,\, t\in \Gamma.
\end{equation}
	
	\begin{lemma}\label{lem:MeanFromSequence}
		Let $\mb F$ be a sequence (or net) of subsets of $\Gamma$.  Then there is a mean $m$ on $\ell^\infty(\Gamma)$ satisfying
		\begin{align}\label{eq:EfmEf}
			\underline{\mathbb E}_{\mb F}(f)\leq &\,m(f)\leq \bar{\mathbb E}(f) && \text{for all real valued } f\in \ell^{\infty}(\Gamma),\\
	\label{eq:LowerUpperMean}
			\underline{d}_{\mb F}\leq & \,m(A)\leq \bar{d}_{\mb F}(A) && \text{for all } A\subseteq \Gamma.
		\end{align}  If $\mb F$ is a F{\o}lner sequence or F{\o}lner net,  there is such an $m\in \mathcal M_{\tau}(\Gamma)$ satisfying (\ref{eq:EfmEf}) and (\ref{eq:LowerUpperMean}).
	\end{lemma}
	
	\begin{proof}		
		Let $\mb F=(F_{j})_{j\in I}$ be a sequence (or net) of finite subsets of $\Gamma$.  Consider the means $\lambda_{j}\in \mathcal M_{\tau}(\Gamma)$, $j\in I$ given by $\lambda_{j}(f):=\frac{1}{|F_{j}|}\sum_{\gamma\in F_{j}}f(\gamma)$.  Let $m$ be any weak$^{*}$ cluster point of $(\lambda_{j})_{j\in I}$, so that $\liminf_{j}(\lambda_{j}(f))\leq m(f)\leq \limsup_{j}(\lambda_{j}(f))$ for all real-valued $f\in \ell^{\infty}(\Gamma)$.  This simplifies to (\ref{eq:EfmEf}), and (\ref{eq:LowerUpperMean}) is the special case where $f=1_{A}$.

		When $\mb F$ is a F{\o}lner sequence or net, it is straightforward to verify that for every $f\in \ell^{\infty}(\Gamma)$ and $\gamma\in \Gamma$, we have $\lim_{j}(\lambda_{j}(f_{\gamma})-\lambda_{j}(f))=0$. It follows that every weak$^{*}$ cluster point of the net (or sequence) $(\lambda_{j})_{j\in I}$ is an invariant mean.  
	\end{proof}

\begin{remark}
Somewhat surprisingly, \cite[Theorem 5.9]{Hopfensperger_WhenIsAMean} shows that every invariant mean on every discrete abelian group (and in fact, every discrete group) can be obtained as $\mathbb E_{\mb F}$ for some F{\o}lner net $\mb F$.  This was previously obtained for the group $\mathbb Z$ in \cite{HindmanStrauss_DensityAndInvariant}.  While this fact is not needed for our results, it may be useful for applications.
\end{remark}

\subsection{Invariant means and almost periodicity}

It is well known that for every uniformly almost periodic function $\phi$ and any two invariant means $m$, $m'$ on $\ell^{\infty}(\Gamma)$, $m(\phi)=m'(\phi)$. This is a consequence of the fact that $m(\chi)=0$ for every invariant mean $m$ and every nontrivial $\chi\in \widehat{\Gamma}$ (see Lemma \ref{lem:RequiredProperties} below). Linearity of $m-m'$ then implies $m(\phi)=m'(\phi)$ for every trigonometric polynomial $\phi$, and boundedness of $m-m'$  implies $m(\phi)=m'(\phi)$ every uniform limit of trigonometric polynomials.

We use the following consequence in the next example: if $\phi:\mathbb Z\to \mathbb C$ is uniformly almost periodic, $(F_{n})_{n\in \mathbb N}$ is a F{\o}lner sequence, and $m$ is any invariant mean on $\ell^{\infty}(\mathbb Z)$, then
\begin{equation}\label{eq:FolnerIndependent}
	\lim_{n\to\infty} \frac{1}{|F_{n}|}\sum_{x\in F_{n}} \phi(x)=m(\phi)
\end{equation}
	
	\begin{observation} If $H\leq \Gamma$ is a finite index subgroup and $C\subseteq \Gamma$ then $1_{C+H}$ is uniformly almost periodic.  In fact it is a trigonometric polynomial: if $\rho:\Gamma \to \Gamma/H$ is the quotient map, then $1_{C+H}=\sum_{\chi\in \widehat{\Gamma/H}} c_{\chi} \chi\circ \rho$ for some $c_{\chi}\in \mathbb C$, and $\chi \circ \rho\in \widehat{\Gamma}$ for each $\chi$.  Since $\Gamma/H$ is finite, the sum is finite, and this means $1_{C+H}$ is a trigonometric polynomial.
\end{observation}	
	
\begin{example}\label{ex:NeedExtreme}
There are sets $A$, $B\subseteq \mathbb Z$ and an invariant mean $m$ on $\ell^{\infty}(\mathbb Z)$ such that $m(A)>0$, $m(B)>0$, and for every Bohr neighborhood $U\subset \mathbb Z$, $m(U\cap (A+B))<m(U)$.  Consequently, the assertion ``$U\subset_{m} A+B$'' is false for every Bohr neighborhood $U\subseteq \mathbb Z$.
\end{example} 
	To construct the example, fix strictly increasing sequences of natural numbers $a_{n}$, $b_{n}$ such that for all $n\in \mathbb N$, we have
	\begin{align}
	\label{eq:anplus1}	a_{n+1}&>2b_{n},\\
\label{eq:2an}		2a_{n}&>b_{n},\\
\label{eq:bnplus1}		(b_{n+1}-a_{n+1})/b_{n}&\to \infty.
	\end{align} 
Let $q_{n}= \lfloor (b_{n}-a_{n})/2 \rfloor$, and set $I_{n}=[a_{n},a_{n}+q_{n}]\cap \mathbb Z$, $J_{n}=[a_{n}+q_{n}+1,b_{n}]\cap \mathbb Z$, and $F_{n}=[a_{n},b_{n}]$, so that $\mb I=(I_{n})_{n\in \mathbb N}$, $\mb J=(J_{n})_{n\in \mathbb N}$, and $\mb F= (F_{n})_{n\in \mathbb N}$ are F{\o}lner sequences.  
	
	Let $A_{n}=(I_{n}\cap 2\mathbb Z)\cup (J_{n}\cap (2\mathbb Z+1))$, and let $B_{n}=F_{n}\cap 2\mathbb Z$.  Let $A=\bigcup_{n\in \mathbb N} A_{n}$, $B=\bigcup_{n\in \mathbb N} B_{n}$.   Conditions (\ref{eq:anplus1})-(\ref{eq:bnplus1}) imply
	\begin{equation}\label{eq:ABFn}\lim_{n\to\infty} |((A+B)\cap F_{n}) \triangle (A_{n}\cap F_{n})|/|F_{n}|=0.
		\end{equation}
	Let $m$ be an invariant mean obtained by applying Lemma \ref{lem:MeanFromSequence} to $\mb F$. Assume, to get a contradiction, that $U$ is a Bohr neighborhood with $U\subset_{m} A+B$.  By Lemma \ref{lem:BohrEquivalents}, there is a uniformly almost periodic function $\phi: \mathbb Z\to [0,1]$ supported on $U$, with $m(\phi)>0$.  This implies $m(\phi 1_{A+B})=m(\phi 1_{U})=m(\phi)>0$.  Note that $\phi1_{2\mathbb Z}$ and $\phi 1_{2\mathbb Z+1}$ are also uniformly almost periodic. We then have
	\begin{align*}
		m(\phi 1_{A+B})&\leq \limsup_{n\to\infty} \frac{1}{|F_{n}|}\sum_{x\in F_{n}} \phi(x)1_{A+B}(x)  && \text{by Lemma \ref{lem:MeanFromSequence}}\\
		&\leq \lim_{n\to\infty} \frac{1}{|F_{n}|}\sum_{x\in A_{n}} \phi(x)  && \text{by (\ref{eq:ABFn})}\\
		&=\lim_{n\to\infty} \frac{1}{2|I_{n}|}\sum_{x\in I_{n}}\phi(x)1_{2\mathbb Z}(x) +  \frac{1}{2|J_{n}|}\sum_{x\in J_{n}}\phi(x)1_{2\mathbb Z+1}(x)  && \text{since }|F_{n}|=2|I_{n}|=2|J_{n}|\\
		&= (1/2)m(\phi1_{2\mathbb Z})+(1/2)m(\phi 1_{2\mathbb Z+1})  && \text{by (\ref{eq:FolnerIndependent})}\\
		&= (1/2)m(\phi),
	\end{align*}
	so $0<m(\phi)\leq  (1/2)m(\phi)$. This is the desired contradiction.

\subsection{Small sumsets and subgroups}

Corollary \ref{cor:FiniteIndex} demonstrates the utility of Theorem \ref{th:Main}: statements about sumsets in compact abelian groups can be transferred to statements about sumsets in discrete abelian groups.

\begin{corollary}\label{cor:FiniteIndex}
	Let $\Gamma$ be a discrete abelian group and let $\nu\in \mathcal M_{\tau}^{ext}(\Gamma)$,  $\eta\in \mathcal M_{\tau}(\Gamma)$,	and 
$A$, $B\subseteq \Gamma$. If  $\nu(A+B)<\nu(A)+\eta(B)$, then there is a finite index subgroup $H\leq \Gamma$ and $t\in \Gamma$ such that $t+H\subset_{\nu} A+B$.
\end{corollary}

\begin{proof}
	Let $\Gamma$, $\nu$, $\eta$, $A$, and $B$ be as in the hypothesis, and let $\tilde{A} =\tilde{A}_{\nu}$, $ \tilde{B}=\tilde{B}_{\eta}\subseteq b\Gamma$ be the sets provided by Theorem \ref{th:Main}. We will prove that $(\tilde{A}+\tilde{B})\cap \Gamma \subset_{\nu} A+B$ and that $(\tilde{A}+\tilde{B})\cap \Gamma$ is a union of cosets of a finite index subgroup of $\Gamma$. First we claim  
	\begin{equation}\label{eq:critical}
		\mu_{b\Gamma}(\tilde{A}+\tilde{B})<\mu_{b\Gamma}(\tilde{A})+\mu_{b\Gamma}(\tilde{B}).
	\end{equation}  To see this, apply Theorem \ref{th:Main} to find
	\begin{align*}\mu_{b\Gamma}(\tilde{A}+\tilde{B})&\leq \nu(A+B) && \\
		&< \nu(A)+\eta(B)  && \text{by hypothesis} \\ &\leq \mu_{b\Gamma}(\tilde{A})+\mu_{b\Gamma}(\tilde{B}).
	\end{align*}
	Because of (\ref{eq:critical}), Satz 1 of \cite{Kneser_SummenmengenLokalkompakten} provides a compact open (hence finite index) subgroup  $\tilde{H}\leq b\Gamma$ satisfying $\tilde{A}+\tilde{B}=\tilde{A}+\tilde{B}+\tilde{H}$.  Thus $\tilde{A}+\tilde{B}$ is compact and is a union of cosets of $\tilde{H}$.  We may take $V=\tilde{A}+\tilde{B}$ in Theorem \ref{th:Main} and conclude that $(\tilde{A}+\tilde{B})\cap \Gamma \subset_{\nu} A+B$. Setting $H=\tilde{H}\cap \Gamma$, Lemma \ref{lem:FiniteIndex} says that $H$ has finite index in $\Gamma$.  The same lemma implies $(\tilde{A}+\tilde{B})\cap \Gamma$ is a union of cosets of $H$, so $t+H\subset_{\nu} A+B$ for some $t\in \Gamma$.
\end{proof}

\subsection{Acknowledgements}  Thanks are due to Ryan Alweiss for discussions that inspired this article, and to Gabe Conant for comments on an earlier draft.

	\section{Proof of the Steinhaus Lemma}\label{sec:SteinhausProof}
As motivation for our proof of Theorem \ref{th:Main} we prove Theorem \ref{th:SteinhausCompactAbelian} below.  

\subsection{Fourier coefficients and convolutions}
As above we draw from Chapter 4 of \cite{Folland_CourseInAHA} and Chapters 1 and 2 of \cite{Rudin_Fourier} for background.	

Let $K$ be a compact abelian group with Haar probability measure $\mu$.   For $f\in L^{2}(\mu)$, the \emph{Fourier transform} $\hat{f}:\widehat{K}\to \mathbb C$ is defined as $\hat{f}(\chi):=\int f \overline{\chi} \, d\mu$; individual values $\hat{f}(\chi)$ are called \emph{Fourier coefficients} of $f$. The characters form an orthonormal subset of $L^2(\mu)$ spanning a dense subspace of $L^2(\mu)$.  The Fourier transform maps $L^{2}(\mu)$ onto $\ell^{2}(\widehat{K})$ isometrically: when $f$, $g\in L^{2}(\mu)$, we have the Parseval identity
\begin{equation}\label{eq:Parseval}
	\sum_{\chi\in\widehat{K}} \hat{f}(\chi)\overline{\hat{g}(\chi)} = \int f\bar{g}\, d\mu.
\end{equation}
Specializing with $g=f$, we have 
\begin{equation}\label{eq:Bessel}\sum_{\chi\in \widehat{K}}|\hat{f}(\chi)|^2=\int |f|^2\, d\mu,
\end{equation} and the \emph{Fourier series} $\sum_{\chi\in\widehat{K}} \hat{f}(\chi)\chi$ converges to $f$ in the norm topology of $L^2(\mu)$.  Here when we say a series $\sum_{j\in I} a_j$ converges to $a$ in some topology, we mean convergence in the sense of unordered sums: for every neighborhood $O$ of $a$, there is a finite set $I_O\subseteq I$ such that $\sum_{j\in J} a_j\in O$ whenever $J$ is a finite set containing $I_O$.

\begin{definition}
	Let $f, g \in L^2(\mu)$.  The \emph{convolution} $f*g:K\to \mathbb C$ is defined by
	\[
	f*g(x):=\int f(t)g(x-t)\, d\mu(t).
	\]
\end{definition}
We have the standard identity \cite[Theorem 1.2.4 (b)]{Rudin_Fourier}
\begin{equation}\label{eq:ConvolutionToMultiplication}
	\widehat{f*g}(\chi)=\hat{f}(\chi)\hat{g}(\chi) \qquad \text{for all } \chi\in\widehat{K}.
\end{equation}
The usual proof of (\ref{eq:ConvolutionExpansion}) uses Fubini's theorem to interchange the order of an iterated integral involving Haar measure $\mu$ on $K$.  In our proof of Theorem \ref{th:Main}, we will have a finitely additive measure in place of Haar measure, so Fubini's theorem will not be available.  In \S\ref{sec:Outline} we explain how using an \emph{extreme} invariant mean recovers what we need from Fubini's theorem.

\subsection{Notation} If $f:D\to \mathbb R$ is a function and $c\in \mathbb R$, $\{f>c\}$ denotes the level set $\{x\in D: f(x)>c\}$.

	\begin{theorem}\label{th:SteinhausCompactAbelian}Let $K$ be a compact abelian group with Haar probability measure $\mu$ and let $A$, $B\subseteq K$ have $\mu(A)>0$ and  $\mu(B)>0$.  Then $A+B$ contains a nonempty open set.
	\end{theorem}

To prove Theorem \ref{th:SteinhausCompactAbelian}, let $A, B\subseteq K$ be $\mu$-measurable sets.  We average translates of $1_B$ by elements of $A$, where ``average'' is with respect to Haar measure:
	\[
	F(x):=\int 1_{t+B}(x)1_A(t)\, d\mu(t).
	\]
	We see that $1_{t+B}(x)1_{A}(t)>0$ only if  $t\in A$ and $x\in t+B$, meaning $x\in A+B$.  Thus $F$ is supported on $A+B$. Simplifying, we get  $F(x)=\int 1_{A}(t)1_B(x-t)\, d\mu(t)=1_{A}*1_{B}(x)$. To prove Theorem \ref{th:SteinhausCompactAbelian}, it therefore suffices to prove the following lemma.  From it we will conclude that  $1_A*1_B$ is continuous and not identically $0$; the level set $\{1_{A}*1_{B}>0\}$ will then be the desired open subset of $A+B$.
	
\begin{lemma}\label{lem:ConvolutionIsContinuous}
	If $f, g\in L^{2}(\mu)$, then $f*g$ is continuous, and $\int f*g\, d\mu=\int f\, d\mu \int g \, d\mu$.
\end{lemma}	
	
The proof uses the identity $\int h(-t)\, d\mu(t)=\int h(t)\, d\mu(t)$, which follows from uniqueness of Haar measure.	
	
\begin{proof}  Note that when $f,g\in L^2(\mu_K)$, (\ref{eq:Bessel}) implies $\hat{f},\hat{g}\in \ell^{2}(K)$.  Then (\ref{eq:ConvolutionToMultiplication}) and Cauchy-Schwarz imply $\sum_{\chi\in\widehat{K}}|\widehat{f*g}(\chi)|$ converges.  Hence the Fourier series $\sum_{\chi\in\widehat{K}} \hat{f}(\chi)\hat{g}(\chi)\chi$ converges uniformly, and its sum is continuous.  To prove $f*g$ is continuous,  we show that $f*g$ is identically equal to the sum of its Fourier series. This follows from (\ref{eq:Parseval}): writing $\tilde{g}_x(t)$ for $\overline{g(x-t)}$, we have  
\begin{align*}
\widehat{\tilde{g}_x}(\chi)&= \int \overline{g(x-t)}\overline{\chi(t)}\, d\mu(t)\\
&= \int \overline{g(-t)} \overline{\chi(x+t)}\, d\mu(t)  && \text{replacing $t$ with $t+x$}\\
&= \int \overline{g(-t)}\overline{\chi(t)} \overline{\chi(x)}\, d\mu(t)  \\
&= \int \overline{g(t)}\overline{\chi(-t)} \, d\mu(t)\, \overline{\chi(x)} && \text{replacing $t$ with $-t$}, \text{ factoring $\overline{\chi(x)}$}\\
&= \overline{\hat{g}(\chi)} \overline{\chi(x)},
\end{align*}
so $\widehat{\tilde{g}_x}(\chi)=\overline{\hat{g}(\chi)} \overline{\chi(x)}$
for all $x\in K$, $\chi\in\widehat{K}$. Then 
\begin{equation}\label{eq:ConvolutionExpansion}
f*g(x) =	\int f(t) g(x-t)\, d\mu(t) = \int f \overline{\tilde{g}_{x}}\, d\mu = \sum_{\chi\in\widehat{K}} \hat{f}(\chi)\overline{\widehat{\tilde{g}_x}(\chi)}=\sum_{\chi\in\widehat{K}} \hat{f}(\chi)\hat{g}(\chi)\chi(x)
\end{equation}
for all $x\in K$.	 The identity $\int f*g\, d\mu = \int f\, d\mu \int g\, d\mu$ is the special case of (\ref{eq:ConvolutionToMultiplication}) where $\chi$ is the trivial character. \end{proof}

\section{Fourier analysis with means}
In this section we state our second main result, Theorem \ref{th:MainGeneral}.   It has a weaker hypothesis than Theorem \ref{th:Main}, which we explain in the next subsection.  While much of this discussion can be viewed as a special case of the setup in \cite{BjorklundFish_ProductSet} and in \cite{Bjorklund_ProductSetMeasured}, which deal with nonabelian and even nonamenable groups, the specialization to discrete abelian groups here admits some simplifications.

\subsection{The mean ergodic theorem}

For a Hilbert space $\mathcal H$ and a unitary operator $U:\mathcal H\to \mathcal H$, we write $\mathcal H_{U\text{-inv}}$ for the closed subspace $\{w\in \mathcal H: Uw=w\}$. We write $P_{U\text{-inv}}$ for orthogonal projection onto $\mathcal H_{U\text{-inv}}$. The mean ergodic theorem is often stated as follows.
\begin{theorem}\label{th:MeanErgodic}
	Let $U:\mathcal H\to \mathcal H$ be a unitary operator and $w\in \mathcal H$.  Then 
\[	\lim_{N\to\infty} \frac{1}{N}\sum_{n=1}^{N} U^{n} w = P_{U\text{-inv}}w\] in the norm topology of $\mathcal H$.
\end{theorem}
Weak convergence is often all that is required for applications: we have
\[
\lim_{N\to\infty} \frac{1}{N}\sum_{n=1}^{N} \langle v,U^{n}w\rangle=\langle v, P_{U\text{-inv}}w\rangle \qquad \text{for all } v, w\in \mathcal H.
\] 
There are many examples of sequences $(a_{n})_{n\in \mathbb N}$ that can be used in place of $a_{n}=n$ while retaining convergence to $P_{U\text{-inv}}w$.  We say that a sequence of integers $(a_{n})_{n\in \mathbb N}$ is an \emph{ergodic sequence}\footnote{But see \S\ref{sec:HartmanUD} below for a list of alternative names for this concept.} if for every unitary operator on a Hilbert space $\mathcal H$ and all $v, w\in \mathcal H$
\[
\lim_{N\to\infty} \frac{1}{N}\sum_{n=1}^{N} \langle v ,U^{a_{n}}w\rangle = \langle v, P_{U\text{-inv}} w\rangle.
\]
Many examples of ergodic sequences are known -- see \S\ref{sec:HartmanUD} for  discussion and references.

\subsection{Matrix coefficients and FS-means}\label{sec:FSmeans} The mean ergodic theorem remains true when generalized to arbitrary (possibly uncountable) abelian groups, with the minor caveat that one must average over F{\o}lner nets, rather than sequences.  We will need only weak convergence, so we use means instead of averaging over nets of finite subsets of $\Gamma$.

\begin{definition}
	A \emph{unitary action} of $\Gamma$ on a Hilbert space $\mathcal H$ is an action of $\Gamma$ on $\mathcal H$ by unitary operators $U_{\gamma}$, $\gamma \in \Gamma.$ 
	For $v, w\in \mathcal H$, the \emph{matrix coefficient} determined by $U$, $v$, and $w$ is the function $\phi_{v,w}:\Gamma\to \mathbb C$, 
	\[\phi_{v,w}(\gamma):=\langle v, U_{\gamma }w\rangle.\] 
	We write $\mathcal{H}_{U\text{-inv}}$ for the closed subspace $\{x\in \mathcal H: U_{\gamma} x = x \text{ for all } \gamma\in \Gamma\}$, and write $P_{U\text{-inv}}w$ for the orthogonal projection of $w$ onto $\mathcal H_{U\text{-inv}}$.
\end{definition}
The following is a weak ergodic theorem for abelian groups.
\begin{theorem}\label{th:WeakErgodic}
Let $\Gamma$ be a discrete abelian group.	If $m\in \mathcal M_{\tau}(\Gamma)$, $U$ is a unitary action of $\Gamma$ on a Hilbert space $\mathcal H$, and $v, w\in \mathcal H$, then 
\begin{equation}
	m(\phi_{v,w})=\langle v, P_{U\text{-inv}} w\rangle.
\end{equation}
\end{theorem} 
Although Theorem \ref{th:WeakErgodic} is classical (and in fact a generalization to \emph{arbitrary} groups goes back at least to \cite{Godement}), we give a proof in \S\ref{sec:Continuous}.  The conclusion of Theorem \ref{th:WeakErgodic} provides all properties of $m$ required for Theorem \ref{th:Main}, motivating the following definition.

\begin{definition}
	Let $\Gamma$ be a discrete abelian group and let $\eta$ be a (not necessarily invariant) mean on $\ell^{\infty}(\Gamma)$.  We say that $\eta$ is an \emph{FS-mean} if for every unitary action $U$ of $\Gamma$ on a Hilbert space $\mathcal H$ and all $v,w\in \mathcal H$, we have $\eta(\phi_{v,w})=\langle v, P_{U\text{-inv}} w\rangle$. 
\end{definition}
FS-means can be thought of as generalizations of ergodic sequences. Every invariant mean is an FS-mean, by Theorem \ref{th:WeakErgodic}.  We write $\mathcal M^{FS}(\Gamma)$ for the set of $FS$-means on $\Gamma$.  So we have
\[
\mathcal M_{\tau}^{ext}(\Gamma)\subseteq \mathcal M_{\tau}(\Gamma) \subseteq \mathcal M^{FS}(\Gamma).
\]
In countable abelian groups there are many non\text{-}invariant FS-means.  See \S\ref{sec:HartmanUD} for details.

\subsection{Eigenspaces}\label{sec:Eigenspaces}
For a unitary action $U$ of $\Gamma$ on a Hilbert space $\mathcal H$ and  $\chi\in \widehat{\Gamma}$, the \emph{$\chi$-eigenspace of $U$} is 
\[
\mathcal H_{\chi}:=\{w\in \mathcal H: U_{\gamma} w=\chi(\gamma) w \text{ for all } \gamma \in \Gamma\}.
\] 
We write $\mathcal H_{c}$ to denote the smallest closed subspace of $\mathcal H$ containing every eigenspace of $\mathcal U$.  In other words, $\mathcal H_{c}$ is the closure of the span of $\bigcup_{\chi\in \widehat{\Gamma}} \mathcal H_{\chi}$.

The next lemma states the properties of FS-means we use in the proof of Theorem \ref{th:Main}.
\begin{lemma}\label{lem:RequiredProperties}
	Let $\eta$ be a mean on $\ell^{\infty}(\Gamma)$. The following are equivalent.	
	\begin{enumerate}
		\item $\eta$ is an $FS$-mean.
		\item Both of the following hold:
		
		\begin{enumerate} \item\label{item:HartmanDef} For every nontrivial $\chi\in\widehat{\Gamma}$, $\eta(\chi)=0$.
		
		\item\label{item:Annihilates} For every unitary action $U$ of $\Gamma$ on a Hilbert space $\mathcal H$ and all $v$, $w\in \mathcal H$ with $w\perp \mathcal H_{c}$, we have   $\eta(|\phi_{v,w}|^{2})=0$.
	\end{enumerate}
	\end{enumerate} 
\end{lemma}

Lemma \ref{lem:RequiredProperties} is proved as part of Lemma \ref{lem:FSequivalents} in \S\ref{sec:Continuous}.

\begin{corollary}\label{cor:Bessel}
	If $\eta\in \mathcal M^{FS}(\Gamma)$ and $f\in \ell^{\infty}(\Gamma)$, then $\sum_{\chi\in\widehat{\Gamma}} |\eta(f\overline{\chi})|^{2}\leq \eta(|f|^{2})$.
\end{corollary}
\begin{proof}
Let $\eta\in \mathcal M^{FS}(\Gamma)$. Note that the map $\langle \cdot, \cdot \rangle_{\eta}:\ell^{\infty}(\Gamma)\times \ell^{\infty}(\Gamma) \to \mathbb C$ given by $\langle f, g \rangle_{\eta}=\eta(f\overline{g})$  is a nonnegative Hermitian form.  Lemma \ref{lem:RequiredProperties} implies $\eta(\chi\overline{\psi})=0$ whenever $\chi\neq \psi\in \widehat{\Gamma}$,  while $\eta(\chi\overline{\chi})=\eta(1_{\Gamma})=1$, so $\widehat{\Gamma}$ forms an orthonormal set with respect to $\langle \cdot, \cdot \rangle_{\eta}$.  Bessel's inequality then implies $\sum_{\chi\in\widehat{\Gamma}} |\langle f,\chi\rangle_{\eta}|^{2}\leq   \eta(|f|^{2})$.
\end{proof}

\begin{remark} The ``FS'' in ``FS-mean'' refers to the \emph{Fourer-Stieltjes algebra} $\mathcal B(\Gamma)$.  This the smallest uniformly closed subalgebra of $\ell^{\infty}(\Gamma)$ containing every matrix coefficient.  It is well known that $\mathcal B(\Gamma)$ is a translation-invariant subalgebra of $WAP(\Gamma)$, the algebra of weakly almost-periodic functions on $\Gamma$.  In every group $\Gamma$ (abelian or otherwise), there is a unique invariant mean $m_{W}$ on $WAP(\Gamma)$, hence a unique invariant mean $m_{\mathcal B}$ on $\mathcal B(\Gamma)$.  See Section 14 of \cite{Eberlein_AbstractErgodic} for a proof in the case where $\Gamma$ is abelian;  \cite[Corollary 1.26]{Burckel_WAP_book} or \cite[Theorem 1.20]{Burckel_WAP_Thesis}  for a proof in the general case.  Hence if $m$ is an invariant mean on $\ell^\infty(\Gamma)$ and $\phi$ is a matrix coefficient, then $m(\phi)=m_{\mathcal B}(\phi)$.
\end{remark}

Lemma \ref{lem:RequiredProperties} is often derived as a consequence of Wiener's lemma for Fourier transforms of measures -- see \cite{Eberlein_FS} for such a statement in the generality of our setting. We present essentially this argument in our proof of Lemma \ref{lem:FSequivalents}.  Lemma \ref{lem:RequiredProperties} also follows from  Corollary 3.9 of \cite{Badea_Grivaux_KazhdanED}, an explicit computation of $m_{\mathcal B}(\phi_{v,w})$. The same computation is done in Lemma 6.4 and Proposition 6.7 of \cite{BjorklundFish_ProductSet}.  Although the latter article makes the standing assumption that the acting group is countable, that assumption is not used for these results.  

Like Theorem \ref{th:WeakErgodic}, the natural generalization of Lemma \ref{lem:RequiredProperties} to arbitrary (including uncountable and nonamenable) groups goes back to \cite{Godement}.

\subsection{Fourier coefficients, \texorpdfstring{$b\Gamma$}{bGamma}, and the Besicovitch-Radon-Nikodym density}

Let $\Gamma$ be a discrete abelian group.

	\begin{definition}\label{def:meanFourierCoefficients} Given $f\in \ell^\infty(\Gamma)$ and a (not necessarily invariant) mean $m$ on $\ell^\infty(\Gamma)$, the \emph{Fourier coefficients of $f$ with respect to $m$} are $\hat{f}^m(\chi):=m(f\overline{\chi})$, where $\chi\in\widehat{\Gamma}$.  
\end{definition}

Recall from \S\ref{sec:BohrCompactification} that we identify $\Gamma$ with a topologically dense subgroup of its Bohr compactification $b\Gamma$, and that every character $\chi\in \widehat{\Gamma}$ is the restriction of a unique $\tilde{\chi}\in \widehat{b\Gamma}$.

\begin{lemma}\label{lem:CorrectAveragingAndbGamma}
	Let $\eta$ be an FS-mean on $\ell^\infty(\Gamma)$ and let $\tilde{\phi}:b\Gamma\to \mathbb C$ be continuous.  Then $\int \tilde{\phi}\, d\mu_{b\Gamma}=\eta(\tilde{\phi}|_\Gamma)$.
\end{lemma}

\begin{proof}
	It suffices to prove the identity when $\tilde{\phi}$ is a trigonometric polynomial, as these form a uniformly dense subset of $C(b\Gamma)$.  Writing  $\tilde{\phi}=\sum_{\chi\in\widehat{\Gamma}} c_\chi \tilde{\chi}$,  this reduces to proving that $\eta(\chi)=0$ for all nontrivial $\chi\in \widehat{\Gamma}$, which follows from Lemma \ref{lem:RequiredProperties}.
\end{proof}

\begin{lemma}\label{lem:DefBRN}
	If $\eta$ is an FS-mean on $\ell^\infty(\Gamma)$ and $f\in \ell^{\infty}(\Gamma)$, then  $\hat{f}^{\eta}\in \ell^{2}(\widehat{\Gamma})$, and 
	\[
	 \tilde{f}^{(\eta)}:=\sum_{\chi\in \widehat{\Gamma}} \hat{f}^\eta(\chi)\tilde{\chi}
	\] 
	defines an element of $L^{2}(\mu_{b\Gamma})$ with the following properties:
	\begin{enumerate}
		\item\label{item:Coefficients} $\widehat{\tilde{f}^{(\eta)}}(\tilde{\chi})=\hat{f}^{\eta}(\chi)$ for all $\chi\in\widehat{\Gamma}$;
		\item\label{item:ContinuousIntegral}  $\int \tilde{\phi}\tilde{f}^{(\eta)}\, d\mu_{b\Gamma} = \eta(\tilde{\phi}|_{\Gamma} f)$ for all $\tilde{\phi}\in C(b\Gamma)$;
		\item\label{item:Integral} $\int \tilde{f}^{(\eta)}\, d\mu_{b\Gamma}=\eta(f)$;
		\item\label{item:BRNbound} $\|\tilde{f}^{(\eta)}\|_{L^{\infty}(\mu_{b\Gamma})}\leq \|f\|_{\infty}$. 		Furthermore, if $f$ is real-valued, then $\tilde{f}^{(\eta)}$ is real-valued and $\inf f \leq \tilde{f}^{(\eta)}(x) \leq \sup f$ for $\mu_{b\Gamma}$-almost every $x$.\\  In particular, if $f:\Gamma \to [0,1]$, then $0\leq \tilde{f} \leq 1$ $\mu_{b\Gamma}$-a.e.
	\end{enumerate}
\end{lemma}
We call $\tilde{f}^{(\eta)}$ the \emph{Besicovitch-Radon-Nikodym density of $f$ with respect to $\eta$}; cf.~Section 4 of \cite{BjorklundGriesmer_BohrTriple}.

	\begin{proof} Let $\eta$ and $f$ be as in the hypothesis.  By
	  Corollary \ref{cor:Bessel}, the series $\sum_{\chi\in\widehat{\Gamma}}|\hat{f}^{\eta}(\chi)|^{2}$ converges, so  $\tilde{f}^{(\eta)}$ is well-defined in $L^{2}(\mu_{b\Gamma})$.	
	  
	  Property \ref{item:Coefficients} follows immediately from the definition of $\tilde{f}^{(\eta)}$. Note that \ref{item:Coefficients} is equivalent to the special case of \ref{item:ContinuousIntegral} where $\tilde{\phi}\in \widehat{b\Gamma}$.  The general case of \ref{item:ContinuousIntegral} now follows by linearity and the fact that the the set of trigonometric polynomials is uniformly dense in $C(b\Gamma)$.  \ref{item:Integral} is the special case of \ref{item:Coefficients} where $\chi\in \widehat{\Gamma}$ is  trivial.

We prove the special case of  \ref{item:BRNbound} where $f:\Gamma \to [0,1]$;  the general case follows by linearity. We first prove that 
\begin{equation}\label{eq:in01}
0\leq \int \tilde{\phi} \tilde{f}^{(\eta)} \, d\mu_{b\Gamma}\leq \|\tilde{\phi}\|_{L^{1}(\mu_{b\Gamma})}
\end{equation} whenever $\tilde{\phi}\in L^{1}(\mu)$ is real-valued and nonnegative.   To see this, note that (\ref{eq:in01}) holds for every continuous $\tilde{\phi}:b\Gamma\to [0,1]$, since
\[
	\int \tilde{\phi} \tilde{f}^{(\eta)} \, d\mu_{b\Gamma}=\eta(\tilde{\phi}|_{\Gamma} f)\leq  \eta(\tilde{\phi}|_{\Gamma})=\int \tilde{\phi}\, d\mu_{b\Gamma} = \|\tilde{\phi}\|_{L^{1}(\mu_{b\Gamma})}.
\] 
The general case of (\ref{eq:in01}) follows, as every nonnegative real-valued $\tilde{\phi}\in L^{1}(\mu_{b\Gamma})$  can be approximated in the $L^{1}$-norm by continuous functions taking nonnegative values.  Thus the map $\tilde{\phi}\mapsto \int \tilde{\phi} \tilde{f}^{(\eta)} \mu_{b\Gamma}$ is a linear functional on $L^{2}(\mu_{b\Gamma})$ with norm at most $1$, mapping nonnegative real valued $\tilde{\phi}$ to nonnegative real numbers. The usual identification of $L^{\infty}(\mu_{b\Gamma})$ with the dual of $L^{1}(\mu_{b\Gamma})$ then implies $\tilde{f}^{(\eta)}$ is equal, $\mu_{b\Gamma}$-almost everywhere, to a function taking values in $[0,1]$.
\end{proof}

\begin{remark}
	When $\eta$ is not an FS-mean, one can construct $\tilde{f}^{(\eta)}$ as a \emph{Radon measure} on $b\Gamma$ rather than a function in $L^\infty(\mu_{b\Gamma})$.  It is unclear what assumptions on $\eta$ will lead to useful properties of this measure.  
\end{remark}

 \begin{lemma}\label{lem:OneWayEstimate}
	Let $A\subseteq \Gamma$, and let $\eta$ be an FS-mean on $\ell^{\infty}(\Gamma)$.  Let $f=1_{A}$, and let $\tilde{A}_{\eta}$ be any $\mu_{b\Gamma}$-measurable subset satisfying 
	$\tilde{A}_{\eta} \sim_{\mu_{b\Gamma}} \{\tilde{f}^{(\eta)}>0\}$.
	
	\begin{enumerate} \item\label{item:Vcompact} If $V\subseteq b\Gamma$ is compact, then $\mu_{b\Gamma}(\tilde{A}_\eta\cap V)\geq \eta(A\cap V)$.
		
		\item\label{item:Vclopen} If $V\subseteq b\Gamma$ is clopen, then $\eta(A\cap V)=0$ if and only if $\mu_{b\Gamma}(\tilde{A}_\eta \cap V)=0$.
	\end{enumerate}
\end{lemma}
\begin{remark}
	Since $\tilde{f}^{(\eta)}$ is an element of $L^{\infty}(\mu_{b\Gamma})$, the condition $\tilde{A}_{\eta}\sim_{\mu_{b\Gamma}} \{\tilde{f}^{(\eta)}>0\}$ means that for every function $\psi$  belonging to the $L^{\infty}(\mu_{b\Gamma})$-equivalence class of $\tilde{f}^{(\eta)}$, we have $\mu_{b\Gamma}(\tilde{A}_{\eta}\triangle \{\psi>0\})=0$.
\end{remark}

\begin{proof}
To prove \ref{item:Vcompact}, assume $V\subseteq b\Gamma$ is compact.	By assumption, $\tilde{A}_{\eta}\sim_{\mu_{b\Gamma}} \{\tilde{f}^{(\eta)}>0\}$. 
	Since $0\leq \tilde{f}^{(\eta)}\leq 1$ $\mu_{b\Gamma}$-a.e., it suffices to prove that 
	\begin{equation}\label{eq:int1Ef} 
		\int 1_{V} \tilde{f}^{(\eta)}\, d\mu_{b\Gamma}\geq \eta(A\cap V).
	\end{equation}  
	Lemma \ref{lem:DefBRN} says that $\int \tilde{\phi} \tilde{f}^{(\eta)}\, d\mu_{b\Gamma}=\eta(\tilde{\phi}|_{\Gamma} 1_A)$ for every continuous $\tilde{\phi}:b\Gamma\to [0,1]$.  Thus, for every continuous $\tilde{\phi}:b\Gamma\to [0,1]$ with $\tilde{\phi}|_{V}\equiv 1$, we have $\int \tilde{\phi}\, \tilde{f}^{(\eta)}\, d\mu_{b\Gamma}\geq \eta(A\cap V)$.  By outer regularity of Haar measure, this implies (\ref{eq:int1Ef}).
	
	To prove \ref{item:Vclopen}, assume $V\subseteq b\Gamma$ is clopen.  Part \ref{item:Vcompact} provides the implication $\mu_{b\Gamma}(\tilde{A}_\eta \cap V)=0$ $\implies$ $\eta(A\cap V)=0$.  To prove the reverse implication, assume $\eta(A\cap V)=0$.  Since $V$ is clopen, $1_V$ is continuous, so Lemma \ref{lem:CorrectAveragingAndbGamma} implies $\int 1_V \tilde{f}^{(\eta)}\, d\mu = \eta(1_V|_\Gamma 1_A)=\eta(A\cap V)=0$.  Since $\tilde{A}_\eta \sim_{\mu_{b\Gamma}} \{\tilde{f}^{(\eta)}>0\}$, this means $\mu_{b\Gamma}(\tilde{A}_\eta \cap V)=0$.
\end{proof}

	\begin{lemma}\label{lem:Main}
	Let $\nu\in \mathcal M_{\tau}^{ext}(\Gamma)$ and $\eta\in \mathcal M^{FS}(\Gamma)$.  Let $A$, $B\subseteq \Gamma$,  let $f=1_A$, $g=1_B$, and let $\tilde{f}^{(\nu)}$, $\tilde{g}^{(\eta)}$ be as in Lemma \ref{lem:DefBRN}.  Let $\tilde{h}=\tilde{f}^{(\nu)}*\tilde{g}^{(\eta)}$ and let $h=\tilde{h}|_{\Gamma}$.  Then the series
	\begin{equation}\label{eq:NuEtaConvolutionExpansion}
		\sum_{\chi\in\widehat{\Gamma}} \hat{f}^{\nu}(\chi)\hat{g}^{\eta}(\chi)\chi
	\end{equation}
converges uniformly to $h$, and
	\begin{enumerate}
				\item\label{item:GeneralHAplusB} $\nu(h1_{A+B})=\nu(h)$;
		\item\label{item:GeneralLevelContained} for all $\delta>0$, $\{h\geq \delta\}\subset_\nu A+B$;
		\item\label{item:GeneralClopenSubset} if $V\subseteq \{\tilde{h}>0\}$ is compact, then $V\cap \Gamma \subset_{\nu} A+B$;
		\item\label{item:GeneralMeasureOfAplusB} $\nu(A+B)\geq \mu_{b\Gamma}(\{\tilde{h}>0\})$.
	\end{enumerate}
 \end{lemma}
 Lemma \ref{lem:Main} is proved in \S\ref{sec:ProofOfMain}.

 The next lemma connects the level set $\{\tilde{h}>0\}$ to a sumset in $b\Gamma$.

 \begin{lemma}\label{lem:LevelSetSumset}
 	Let $K$ be a compact abelian group with Haar probability measure $\mu$, and let $\tilde{f}$, $\tilde{g}:K\to [0,1]$ be $\mu$-measurable functions. Then there are $\mu$-measurable sets $\tilde{A}$, $\tilde{B}\subseteq K$ such that $\mu(\tilde{A})\geq \mu(\{\tilde{f}>0\})$, $\mu(\tilde{B})\geq \mu(\{\tilde{g}>0\})$, while  
 	$\tilde{f}(a)>0$, $\tilde{g}(b)>0$, and $\tilde{f}*\tilde{g}(a+b)>0$ for all $a\in \tilde{A}$, $b\in \tilde{B}$.
 	
 	Consequently, $\{\tilde{f}*\tilde{g}>0\}$ is an open neighborhood of $\tilde{A}+\tilde{B}$.
 \end{lemma}
 
 Lemma \ref{lem:LevelSetSumset} is proved in \S\ref{sec:LevelSumset}.  
 
  The next lemma deduces property  \ref{item:GeneralLevelContained} in Lemma \ref{lem:Main} from \ref{item:GeneralHAplusB}.
 
 \begin{lemma}\label{lem:EssentialLevelContainment}
 	Let $m$ be a mean on $\ell^{\infty}(\Gamma)$ and $D\subseteq \Gamma$.  Assume $h:\Gamma \to [0,1]$ satisfies $m(h\phi)\geq 0$ for all $\phi:\Gamma\to[0,1]$, and $m(h1_{D})=m(h)$. Then for all $\delta>0$, $\{h\geq \delta\}\subset_{m} D$.
 \end{lemma}
 
 \begin{proof}
 	Assume $m$, $h$, and $D$ are as in the hypothesis. Let $C:=\Gamma \setminus D$. Assume, to get a contradiction, that $\delta>0$ and $\alpha:=m(1_{C\cap \{h\geq\delta\}})>0$.  Then 
 	\begin{align*}
 		m(h)&= m(h1_{C\cap \{h< \delta\}}) + m(h1_{C\cap \{h\geq \delta\}}) + m(h1_{D})\\ &\geq 0+\delta\alpha+m(h1_{D})\\
 		&>m(h1_{D}),
 	\end{align*} which contradicts the assumption $m(h)=m(h1_{D})$.  
 \end{proof}

 Combining Lemmas \ref{lem:Main} and \ref{lem:LevelSetSumset} yields our main result, the promised generalization of Theorem \ref{th:Main}.
	
	\begin{theorem}\label{th:MainGeneral}  Let $\Gamma$ be a discrete abelian group.  Let $\nu\in \mathcal M_{\tau}^{ext}(\Gamma)$, $\eta\in \mathcal M^{FS}(\Gamma)$, and $A$, $B\subseteq \Gamma$.  
	Then there are $F_{\sigma}$ sets $\tilde{A}_{\nu}$, $\tilde{B}_{\eta}\subseteq b\Gamma$  such that $\mu_{b\Gamma}(\tilde{A}_{\nu})\geq \nu(A)$, $\mu_{b\Gamma}(\tilde{B}_{\eta})\geq \eta(B)$, and
\begin{equation}\label{eq:MainMeasure}
		\nu(A+B)\geq \mu_{b\Gamma}(\tilde{A}_{\nu}+\tilde{B}_{\eta}).
\end{equation}
Furthermore,
\begin{enumerate}
	\item\label{item:recallAVBV} if $V\subseteq b\Gamma$ is compact, then 
	\[\mu_{b\Gamma}(\tilde{A}_{\nu}\cap V)\geq \nu(A\cap V) \quad \text{and} \quad \mu_{b\Gamma}(\tilde{B}_{m}\cap V)\geq m(B\cap V);\]
	
	\item\label{item:recallVbGammaClopen} if $V\subseteq b\Gamma$ is clopen, then
	\[\mu_{b\Gamma}(\tilde{A}_{\nu}\cap V)=0\iff \nu(A\cap V)=0 \quad \text{and} \quad \mu_{b\Gamma}(\tilde{B}_{m}\cap V)=0\iff m(B\cap V)=0;\]
	
	\item\label{item:recallVcapGamma}  if $V\subseteq \tilde{A}_{\nu}+\tilde{B}_{m}$ is compact, then $V\cap \Gamma \subset_{\nu} A+B$.
\end{enumerate}
\end{theorem}

\begin{proof}
	Let $\Gamma$, $\nu$, $\eta$, $A$, and $B$ be as in the hypothesis. Let $\tilde{f}^{(\nu)}$, $\tilde{g}^{(\eta)}$, and $\tilde{h}$ be as in Lemma \ref{lem:Main}. Then $\tilde{f}^{(\nu)}$, $ \tilde{g}^{(\eta)}:b\Gamma\to [0,1]$, so Lemma \ref{lem:LevelSetSumset} provides $\mu_{b\Gamma}$-measurable sets $\tilde{A}_{\nu}\subseteq \{\tilde{f}^{(\nu)}>0\}$ and $\tilde{B}_{\eta}\subseteq \{\tilde{f}^{(\eta)}>0\}$ satisfying $\mu_{b\Gamma}(\tilde{A}_{\nu})\geq \mu_{b\Gamma}(\{\tilde{f}>0\})$, $\mu_{b\Gamma}(\tilde{B}_{\eta})\geq \mu_{b\Gamma}(\{\tilde{g}>0\})$, and $\tilde{A}_{\nu}+\tilde{B}_{\eta}\subseteq \{\tilde{h}>0\}$.   Inner regularity of Haar measure allows us to replace $\tilde{A}_{\nu}$ and $\tilde{B}_{\eta}$ with countable unions of compact sets having the same $\mu_{b\Gamma}$-measure, so that $\tilde{A}_{\nu}$ and $\tilde{B}_{\eta}$ are $F_{\sigma}$ sets.
	
	Inequality (\ref{eq:MainMeasure}) follows from Lemma \ref{lem:LevelSetSumset} \ref{item:GeneralMeasureOfAplusB} and the containment $\tilde{A}_{\nu}+\tilde{B}_{\eta}\subseteq \{\tilde{h}>0\}$. 	Parts \ref{item:recallAVBV} and \ref{item:recallVbGammaClopen} follow from Lemma \ref{lem:OneWayEstimate}.
		
	To prove part \ref{item:recallVcapGamma}, assume $V\subseteq \tilde{A}_{\nu}+\tilde{B}_{\eta}$ is compact. Then $V$ is a compact subset of  $\{\tilde{h}>0\}$. Continuity of $\tilde{h}$ and compactness of $V$ then imply $V\subseteq \{\tilde{h}\geq \delta\}$ for some $\delta>0$.  Part \ref{item:GeneralClopenSubset} of Lemma \ref{lem:Main} then implies $V\subset_{\nu} A+B$. \end{proof}

\begin{proof}[Proof of Theorem \ref{th:Main}]
	Since every invariant mean is an FS-mean, Theorem \ref{th:Main} is a special case of Theorem \ref{th:MainGeneral}.
\end{proof}

	\section{Piecewise Bohr structure in \texorpdfstring{$A+B$}{A+B}}

	Corollary \ref{cor:Recover} generalizes the results discussed in \S\ref{sec:JinBFW} and recovers \cite[Theorem 1.4]{Griesmer_SumsetsDenseSparse}.  To prove it, we need the following general fact about translation\text{-}invariant measures on groups.
	
	\begin{lemma}\label{lem:probability}
		Let $G$ be an abelian  group and let $m$ be a translation invariant, finitely additive (or countably additive) probability measure on $G$ such that $m(-A)=m(A)$ for all $m$-measurable $A\subseteq G$.  If $m(A) + m(B)>1$, then $A+B=G$.
	\end{lemma}
	
	\begin{proof}
		Let $x\in G$.  Note that $x\in A+B$ if and only if $B\cap (x-A)\neq \varnothing$.  Since $m$ is translation invariant and $m(-A)=m(A)$, we have $m(B)+m(x-A)>1$.  Hence $B$ cannot be disjoint from $x-A$, so $x\in A+B$.
	\end{proof}
	
	Note that Lemma \ref{lem:probability} applies to Haar probability measure $\mu$ on a compact abelian group: uniqueness of Haar measure implies $\mu(A)=\mu(-A)$ for every $\mu$-measurable set $A$.
	
	\begin{corollary}\label{cor:Recover}
		Let $\Gamma$ be a discrete abelian group and let $A$, $B\subseteq \Gamma$.
		\begin{enumerate}			
			\item\label{item:RecoverG} If $d^*(A)>0$ and $\eta(B)>0$ for some FS-mean $\eta$ on $\ell^\infty(\Gamma)$, then $A+B$ is piecewise Bohr (and  therefore piecewise syndetic).  Furthermore, if $\nu\in \mathcal M_{\tau}^{ext}(\Gamma)$ and $\nu(A)>0$, then there is a Bohr neighborhood $U\subseteq \Gamma$ such that $U \subseteq_{\nu} A+B$.
	
			\item\label{item:RecoverThick} If $\nu$ is an extreme invariant mean on $\ell^{\infty}(\Gamma)$, $\eta$ is an FS-mean on $\ell^\infty(\Gamma)$, and $\nu(A)+\eta(B)>1$, then $\nu(A+B)=1$. Consequently $A+B$ is thick. 			
\end{enumerate}
	\end{corollary}

	\begin{proof}
		The proof of  \ref{item:RecoverG} is identical to the proof of Corollary \ref{cor:PWB}, with Theorem \ref{th:MainGeneral} in place of Theorem \ref{th:Main}. 
		
		To prove \ref{item:RecoverThick}, assume $\nu$, $\eta$, $A$, and $B$ are as in the hypothesis, and let  $\tilde{A}_{\nu}, \tilde{B}_{\eta}\subseteq b\Gamma$ be as in Theorem \ref{th:MainGeneral}.  Then  $\mu_{b\Gamma}(\tilde{A}_{\nu})\geq \nu(A)$ and $\mu_{b\Gamma}(\tilde{B}_{\eta})\geq \eta(B)$, so $\mu_{b\Gamma}(\tilde{A}_{\nu})+\mu_{b\Gamma}(\tilde{B}_{\eta})> 1$.  Lemma \ref{lem:probability} implies $\tilde{A}_{\nu}+\tilde{B}_{\eta}=b\Gamma$, so that $\mu_{b\Gamma}(\tilde{A}_{\nu}+\tilde{B}_{\eta})=1$. Theorem \ref{th:MainGeneral}  then implies $\nu(A+B)=1$.			Thickness of $A+B$ then follows from  \ref{item:MeanEquals1}$\implies$\ref{item:isThick} in Lemma \ref{lem:ThickEquivalents}. \end{proof}

	\subsection{Connection to previous work}	Lemma \ref{lem:Main} and	Theorem \ref{th:MainGeneral} improve \cite[Theorem I]{BergelsonFurstenbergWeiss}, which says that if $A, B\subseteq \mathbb Z$ have positive upper Banach density, then $A+B$ is piecewise Bohr.  The improvement is in three aspects:
		
		\begin{enumerate}
			\item The containments $\{h\geq \delta\} \subset_\nu A+B$ and $(\tilde{A}_{\nu}+\tilde{B}_{\eta})\cap V \subset_{\nu} A+B $ in Lemma \ref{lem:Main} and Theorem \ref{th:MainGeneral} say that the piecewise Bohr structure of $A+B$ can be found where the largeness of $A$ itself is realized;  \cite[Theorem I]{BergelsonFurstenbergWeiss} does not specify the location of the thick set $T$ where the piecewise Bohr structure of $A+B$ is realized.
			
			\item The hypothesis on $B$ is weakened from ``positive upper Banach density'' to $\eta(B)>0$ for some FS-mean $\eta$.
			
			\item Relating $A+B$ to  $\tilde{A}_{\nu}+\tilde{B}_{\eta}$ allows one to transfer results about sumsets in compact abelian groups directly to discrete abelian groups.
		\end{enumerate}

	Theorem \ref{th:MainGeneral} is analogous to a special case of \cite[Proposition 1.10]{BjorklundFish_ProductSet}.  The latter does specify the structure of $A+B$ in terms of a convolution $\tilde{f}*\tilde{g}$ on $b\Gamma$, and the results there apply to general countable groups (not necessarily abelian, or even amenable).   Our Theorem \ref{th:MainGeneral} locates the piecweise Bohr structure of $A+B$ in terms of $\nu$, whereas \cite[Proposition 1.10]{BjorklundFish_ProductSet} locates the structure of $A+B$ in terms of a thick set $T$ whose relation to $\nu$ is unspecified.   The proofs in \cite{BjorklundFish_ProductSet} rely on detailed analysis of orbits in measure preserving dynamical systems, which our proofs avoid.
	
	Results in \cite{BjorklundFish_ApproxInvariance} do locate structure in $A+B$ in terms of a $\nu\in \mathcal M_{\tau}^{ext}(A)$ with $\nu(A)>0$.
	
	Our results apply to discrete abelian groups of arbitrary cardinality, whereas \cite{BjorklundFish_ProductSet} and \cite{BjorklundFish_ApproxInvariance} consider only countable discrete groups.

	\section{Convolution on discrete groups}\label{sec:Outline}
	
		\subsection{Integral notation for means}  While means are only finitely additive, we follow \cite{BjorklundFish_ProductSet} and \cite{BjorklundFish_ApproxInvariance}, using integral notation to denote some evaluations of means.  So $m(f)$ may be denoted $\int f(t)\, dm(t)$.
	 
\subsection{Pointwise convolution versus weak convolution}	To prove Lemma \ref{lem:Main}, one would naturally set $f=1_A$, $g=1_B$, and define a convolution $f*_\eta g\in \ell^\infty(\Gamma)$ by $f*_\eta g(x):= \int f(x-t)g(t)\, d\eta(t)$.  One easily verifies that $f*_\eta g(x)=0$ unless $x\in A+B$, so we would like to show that
\begin{equation}\label{eq:Want}
	f*_\eta g \sim_\nu \sum_{\chi\in\widehat{\Gamma}} \hat{f}^\nu(\chi) \hat{g}^\eta(\chi) \chi,
\end{equation} cf.~equation (\ref{eq:ConvolutionExpansion}) above. This is difficult,  due the fact that $\nu$ is not countably additive: the inability to exchange the order of integration \`a la Fubini is an obstruction to computing Fourier coefficients of $f*_\eta g$.  To circumvent this issue, we use the alternative Definition \ref{def:nuEtaConvolution}, denoted $f\leftindex_{\nu}*_\eta\, g$. This can be thought of as a \emph{weak $L^2(\nu)$ limit} of averages of translates of $f$, weighted by $g$ and $\eta$.  When $\nu$ is only assumed to be invariant, we still cannot prove an identity  like (\ref{eq:Want}) for $f\leftindex_{\nu}*_\eta\, g$, due to examples like Example \ref{ex:NeedExtreme}.  Thanks to Lemma \ref{lem:TranslationExtremeErgodic},  the desired identity is available if we assume $\nu\in \mathcal M_{\tau}^{ext}(\Gamma)$ and $\eta\in \mathcal M^{FS}(\Gamma)$.
	
	\subsection{Averaging unitary actions} \label{sec:DefE}  We generalize the mean ergodic theorem (Theorem \ref{th:MeanErgodic} above) in two ways: (i) instead of averaging over $n\in \{1,\dots,N\}$, we take averages with respect to a mean on $\ell^\infty(\Gamma)$, and (ii) we introduce some coefficients $g(n)$, so our averages will be generalizations of $\frac{1}{N}\sum_{n=1}^N g(n)U^n v$.  We will see that limits of such averages can often be computed in terms of Fourier coefficients and eigenspaces; cf.~\cite{JonesRosenblattTempelman_Ergodic}.
	
	 Let $\eta$ be any mean on $\ell^{\infty}(\Gamma)$ and  $U$ an action of $\Gamma$ by unitary operators $U_t$ on a Hilbert space $\mathcal H$. For each $w\in \mathcal H$ and $g\in \ell^{\infty}(\eta)$ we define \[\int (U_{t} w)g(t)\, d\eta(t)\] to be the unique $z\in \mathcal H$ satisfying $\langle v, z\rangle = \int \langle v, U_t w\rangle\, g(t) d\eta(t)$ for all $v\in \mathcal H$. The Reisz representation theorem for Hilbert spaces provides this $z$, as the map $v\mapsto \int \langle v, U_t w\rangle\, g(t) d\eta(t)$ is a continuous linear functional on $\mathcal H$.
	
	Recall from \S\ref{sec:Eigenspaces} that for a given $\chi\in \widehat{\Gamma}$, $\mathcal H_\chi:=\{w\in \mathcal H:U_{\gamma}w= \chi(\gamma)w \text{ for all } \gamma\in \Gamma\}$ is the $\chi$-eigenspace of $\mathcal U$.
For each $g\in \ell^{\infty}(\Gamma)$ and each $w\in \mathcal H_{\chi}$, we have $\int (U_{t} w) g(t)\, d\eta(t)=\int  \chi(t)g(t)w \, d\eta(t)$, so 
	\begin{equation}\label{eq:EigenMean}\int (U_{t}w)g(t)\,d\eta(t) = \eta(g\chi)w \qquad \text{ for each } w\in \mathcal H_{\chi}, g\in \ell^{\infty}(\Gamma).
	\end{equation}

	\subsection{Hilbert space associated to a mean} Let $m$ be a mean on $\ell^{\infty}(\Gamma)$. As above we write $f\sim_m g$ if $m(|f-g|)=0$.  
	
	Consider the Hermitian form $\langle \cdot , \cdot \rangle_m$ on $\ell^\infty(\Gamma)$ given by $\langle f,g\rangle_m:=m(f\overline{g})$ and the associated seminorm $\|f\|_{m}:=\langle f,f\rangle_m$.  Observe that $\langle \cdot , \cdot \rangle_m$ is not an inner product, as $\|f\|_m$ may be zero for nonzero $f$.  Nevertheless, $\langle \cdot, \cdot\rangle_m$ is positive semidefinite and sesquilinear; in particular Cauchy-Schwarz holds:
	\begin{equation}\label{eq:CS}
		m(f\overline{g})\leq \|f\|_{m}\|g\|_{m}.
	\end{equation} Let $\ell^{\infty}(\Gamma)/\sim_{m}$ denote the quotient $\ell^\infty(\Gamma)/\{f:\|f\|_m=0\}$.  We write $L^2(m)$ for the metric completion of $\ell^\infty(\Gamma)/\sim_{m}$ with respect to the metric $d_m([f],[g]):=\|f-g\|_m$, where $[f]$ and $[g]$ denote the respective $\sim_{m}$ equivalence classes of $f$ and $g$. 

Elements of $L^2(m)$ are represented by Cauchy sequences $(f_n)_{n\in \mathbb N}$ of elements of $\ell^\infty(\Gamma)$, where ``Cauchy'' means  $\lim_{n\to\infty} \sup_{j,k\geq n} \|f_j-f_k\|_m=0$.  When appropriate, we identify $f\in \ell^{\infty}(\Gamma)$ with the element of $L^{2}(m)$ determined by the constant Cauchy sequence $f_{n}=f$.  In particular, each $\chi\in\widehat{\Gamma}$ determines an element of $L^{2}(m)$.

Now $L^{2}(m)$ is a Hilbert space; the inner product of two elements $v$, $w$ represented by Cauchy sequences $(f_{n})_{n\in \mathbb N}$, $(g_{n})_{n\in \mathbb N}$ is $\langle v,w \rangle_{m}=\lim_{n\to\infty} \langle f_{n},g_{n}\rangle_{m}.$
	
If $w\in L^{2}(m)$ and $g\in \ell^{\infty}(m)$, we write $w\sim_{m} g$ to mean that $w$ is represented by the constant Cauchy sequence $g_{n}=g$.
	
\begin{lemma}\label{lem:FSisHartmanUD}  If $\eta$ is an FS-mean on $\ell^{\infty}(\Gamma)$, then  $\widehat{\Gamma}$ forms an orthonormal subset of $L^{2}(\eta)$.
\end{lemma}

\begin{proof}
This follows from the proof of Corollary \ref{cor:Bessel}.
\end{proof}
	
	\subsection{Translation action}  Now assume $m$ is an invariant mean on $\ell^{\infty}(\Gamma)$.  For $f\in \ell^\infty(\Gamma)$, we write $\gamma\cdot f$ for $f_\gamma$, where  $f_\gamma(x)=f(x-\gamma)$.  Since $m$ is invariant, translation preserves the seminorm: $\|\gamma \cdot f\|_{m}=\|f\|_{m}$ for all $\gamma\in \Gamma$ and $f\in \ell^{\infty}(m)$.   Translation on $\ell^{\infty}(\Gamma)$ therefore induces an action of $\Gamma$ by unitary operators $U_{\gamma}$ on $L^2(m)$: if $w\in L^2(m)$ is represented by the Cauchy sequence $(f_n)_{n\in\mathbb N}$, we define $U_\gamma w$ to be the element of $L^2(m)$ represented by the Cauchy sequence $(\gamma \cdot f_n)_{n\in \mathbb N}$.  It is routine to check that the element of $L^{2}(m)$ determined by $(\gamma \cdot f_n)_{n\in \mathbb N}$ depends only on $w$ and not on the specific choice of Cauchy sequence, so $U_{\gamma}w$ is well-defined.

	\begin{definition}\label{def:nuEtaConvolution}
Let $m$ be an invariant mean and let $\eta$ be any mean on $\ell^{\infty}(\Gamma)$. Let $f,g\in \ell^\infty(\Gamma)$, and write $U$ for the translation action of $\Gamma$ on $L^{2}(m)$. We define $f\leftindex_{m}{*}_\eta\, g$  (the \emph{convolution of $f$ and $g$ with respect to $m$ and $\eta$}) to be the element of $L^{2}(m)$ given by
	\begin{equation}\label{eq:defConvolution}	
	f\leftindex_{m}{*}_\eta\, g := \int (U_t f) g(t)\, d\eta(t).
	\end{equation}
\end{definition}
The ``$f$'' on the right-hand side of (\ref{eq:defConvolution}) is  the element of $L^{2}(m)$ represented by the constant Cauchy sequence $f_{n}=f$, while ``$g$'' denotes the actual function $g$.  Thus $f\leftindex_{m}{*}_\eta\, g$ is the unique element of $L^{2}(m)$ satisfying $ \langle \phi,f\leftindex_{m}{*}_\eta\, g\rangle_{m} = \int \langle \phi, t \cdot f\rangle_{m}\, g(t)\, d\eta(t)$ for all $\phi\in L^{2}(m)$.
	
\begin{lemma}\label{lem:WeakSupportOnAplusB} Let $m$ be an invariant mean on $\ell^\infty(\Gamma)$ and let $\eta$ be any mean on $\ell^\infty(\Gamma)$. If $f,g:\Gamma\to [0,1]$ are supported on $A$ and  $B$, respectively, then $f\leftindex_m*_\eta\, g$ is nonnegative and is supported on $A+B$.  More precisely: if $\phi:\Gamma\to [0,1]$, then
\begin{align}\label{eq:ConvolutionIsPositive}
 \langle \phi, f\leftindex_m*_\eta\, g\rangle_{m} &\geq 0, \\
\label{eq:ActuallyUseful}
		\langle \phi, f\leftindex_m*_\eta\, g \rangle_m &=\langle \phi 1_{A+B}, f\leftindex_m*_\eta\, g\rangle_m .
	\end{align}
\end{lemma}

\begin{proof}
Let $f, g, \phi:\Gamma\to [0,1]$ with $f$ supported on $A$ and $g$ supported on $B$.  For all $b\in B$, we have  $\langle \phi , U_{b}f\rangle_{m} = \langle \phi 1_{A+B} , U_{b}f\rangle_{m}$. Since $g$ is supported on $B$, we get $\int \langle \phi,U_t f\rangle_{m}\, g(t) \, d\eta(t)= \int \langle \phi1_{A+B},U_t f\rangle_{m}\, g(t) \, d\eta(t)$. This implies (\ref{eq:ActuallyUseful}).
	To prove (\ref{eq:ConvolutionIsPositive}), note that  $\langle \phi, U_{t}f\rangle_{m}\geq 0$ for all $t$, so $\int \langle \phi, U_{t}f\rangle_{m}\, g(t)\, d\eta(t)\geq 0$.
\end{proof}

In light of Equation (\ref{eq:EigenMean}), simplifying $f\leftindex_{m}{*}_\eta\, g$ in terms of characters requires us to identify the $\chi$-eigenspace $\mathcal H_{\chi}$ of the translation action $U$ on $L^{2}(m)$.  Clearly we have $\overline{\chi}\in \mathcal H_{\chi}$, but in general $\mathcal H_{\chi}$ may be infinite-dimensional.  Lemma \ref{lem:TranslationExtremeErgodic} says that for \emph{extreme} invariant means $\nu$, each eigenspace of the translation action on $L^{2}(\nu)$ is one-dimensional.
\begin{lemma}\label{lem:TranslationExtremeErgodic}
	Let $\nu\in \mathcal M_{\tau}^{ext}(\Gamma)$ and let $U$ be the translation action of $\Gamma$ on $L^{2}(\nu)$.  Then 
	
	\begin{enumerate} 
		\item\label{item:TranslationErgodicity} the space $\{w\in L^{2}(\nu): U_{\gamma} w = w \text{ for all } \gamma\in \Gamma \}$ of $U$\text{-}invariant elements is one-dimensional and is spanned by $1_{\Gamma}$;
		\item\label{item:TranslationEigenspace} for each $\chi\in \widehat{\Gamma}$, the $\chi$-eigenspace of $U$ is one-dimensional, and is spanned by $\overline{\chi}$.
	\end{enumerate}
\end{lemma}
Lemma \ref{lem:TranslationExtremeErgodic} follows from the well-known correspondence between extreme invariant means on $\ell^{\infty}(\Gamma)$ and Radon probability measures $\mu$ on $\beta\Gamma$ ($=$ the Stone-Cech compactification of the discrete group $\Gamma$)  for which left translation by elements of $\Gamma$ induces an ergodic measure-preserving action on $(\beta\Gamma,\mu)$.  We present an elementary proof in \S\ref{sec:Ergodicity}.

	\begin{lemma}\label{lem:FStarGExpansion}
	Let $\nu\in \mathcal M_{\tau}^{ext}(\Gamma)$ and $\eta\in \mathcal M^{FS}(\Gamma)$. Let $f, g\in \ell^{\infty}(\Gamma)$.  Then the series $\sum_{\chi\in\widehat{\Gamma}} \hat{f}^\nu(\chi)\hat{g}^\eta(\chi)\chi$ converges uniformly to an element $h\in \ell^{\infty}(\Gamma)$ satisfying $f\leftindex_\nu*_\eta\, g \sim_{\nu} h$.
\end{lemma}

Lemma \ref{lem:FStarGExpansion} is proved at the end of \S\ref{sec:SimplifyConvolution}.

\subsection{Simplifying convolutions}\label{sec:SimplifyConvolution} Let $m$ be an invariant mean on $\ell^{\infty}(\Gamma)$ and let $\eta$ be an arbitrary mean on $\ell^{\infty}(\Gamma)$.  Let $U$ be a unitary action of $\Gamma$  on a Hilbert space $\mathcal H$.   

Recall that $\mathcal H_c$ is the smallest closed subspace of $\mathcal H$ containing all the eigenspaces $\mathcal H_{\chi}$.  Since each $U_\gamma$ is unitary, the eigenspaces are closed and mutually orthogonal, and each $w\in \mathcal H$ may be expanded uniquely as 
	\[w=w_0+\sum_{\chi\in\widehat{\Gamma}} w_\chi,\] where $w_\chi\in \mathcal H_\chi$ and $w_0\in \mathcal H_c^\perp$.  With this expansion, equation (\ref{eq:EigenMean}) implies
	\begin{equation}\label{eq:EUetag}
	 \int (U_{t}w) g(t)\, d\eta(t) = 	 \int (U_{t}w_{0}) g(t)\, d\eta(t)+ \sum_{\chi\in\widehat{\Gamma}} \eta(g\chi) w_{\chi}.
	\end{equation}
We will simplify $\int (U_{t}w_{0}) g(t)\, d\eta(t)$ under the assumption that $\eta$ is an FS-mean, and we will simplify $\eta(g\chi)w_{\chi}$ under the assumption that $U$ is the translation action associated to an \emph{extreme} invariant mean.  While these assumptions are independent of one another,  our proof of Lemma \ref{lem:FStarGExpansion} uses both.

\begin{lemma}\label{lem:FSannihilates}
	If $\eta$ is an FS-mean on $\ell^{\infty}(\Gamma)$, $U$ is a unitary action of $\Gamma$ on a Hilbert space $\mathcal H$,  $g\in L^{\infty}(\Gamma)$, and $w\in \mathcal H_{c}^{\perp}$, then $\int (U_{t}w)\, g(t)\, d\eta(t)=0$.
\end{lemma}
\begin{proof}
Let $v\in \mathcal H$ and $w\in \mathcal H_{c}^{\perp}$. Since $\eta$ is an FS-mean, Lemma \ref{lem:RequiredProperties} \ref{item:Annihilates} says that $\int |\langle v, U_{t}w\rangle|^{2}\, d\eta(t)=0$.
We apply Cauchy-Schwarz to find
	\[
	\Bigl|\int \langle v,U_{t}w \rangle g(t) \, d\eta(t)\Bigr| \leq  \int |\langle v, U_{t}w\rangle|^{2}\, d\eta(t)\, \eta(|g|^{2}) =0,
	\]
	so $\int \langle v,U_{t}w \rangle g(t) \, d\eta(t)=0$.  Since this holds for all $v\in \mathcal H$, we have $\int (U_{t}w)\, g(t)\, d\eta(t)=0$.	
\end{proof}

\begin{proof}[Proof of Lemma \ref{lem:FStarGExpansion}]
Let $\nu\in \mathcal M_{\tau}^{ext}(\Gamma)$ and $\eta\in \mathcal M^{FS}(\Gamma)$.  Let $U$ be the translation action on $L^{2}(\nu)$, and let $f, g\in \ell^{\infty}(\Gamma)$.  Specializing equation (\ref{eq:EUetag}) with $f$ in place of $w$, the left-hand side becomes $f\leftindex_{\nu}{*}_{\eta}\, g$, so we have
\begin{equation}\label{eq:FirstExpansion}
	f\leftindex_{\nu}{*}_{\eta}\, g \sim_{\nu} \int (U_{t}f_{0}) g(t)\, d\eta(t) + \sum_{\chi\in \widehat{\Gamma}} \eta(g\chi)f_{\chi},
\end{equation}
where $f_{\chi}$ is the orthogonal projection of $f$ onto the $\chi$-eigenspace of $U$, and $f_{0}$ is orthogonal to every eigenspace of $U$ in $L^{2}(\nu)$.  Lemma \ref{lem:FSannihilates} implies  $\int (U_{t}f_{0}) g(t)\, d\eta(t)=0$. Part \ref{item:TranslationEigenspace} of Lemma \ref{lem:TranslationExtremeErgodic}  implies $f_{\chi}\sim_{\nu}\langle f,\overline{\chi}\rangle_{\nu}\overline{\chi}$, which we can write as $f_{\chi}\sim_{\nu}\nu(f\chi)\overline{\chi}$.  The right-hand side of  (\ref{eq:FirstExpansion}) then simplifies as $\sum_{\chi\in\widehat{\Gamma}} \eta(g\chi)\nu(f\chi)\overline{\chi}$.  Replacing $\chi$ with $\overline{\chi}$ and permuting the order of multiplication, we have $f\leftindex_{\nu}{*}_{\eta}\, g \sim_{\nu}\sum_{\chi\in\widehat{\Gamma}} \nu(f\overline{\chi})\eta(g\overline{\chi})\chi$.  By definition, this is equal to $\sum_{\chi\in\widehat{\Gamma}} \hat{f}^\nu(\chi)\hat{g}^\eta(\chi)\chi$ as desired.

To prove uniform convergence, first apply Corollary \ref{cor:Bessel} to get $\sum_{\chi\in\widehat{\Gamma}}|\eta(g\overline{\chi})|^{2}\leq \eta(|g|^{2})$, and likewise $\sum_{\chi\in\widehat{\Gamma}} |\nu(f\overline{\chi})|^{2}\leq \nu(|f|^{2})$.   Cauchy-Schwarz then implies $\sum_{\chi\in\widehat{\Gamma}} |\nu(f\overline{\chi})\eta(g\overline{\chi})|$ converges, so $\sum_{\chi\in\widehat{\Gamma}} \eta(g\overline{\chi})\nu(f\overline{\chi})\chi$ converges uniformly.
\end{proof}

\section{Ergodicity}\label{sec:Ergodicity}
For this section, fix a discrete abelian group $\Gamma$ and an extreme invariant mean $\nu$ on $\ell^{\infty}(\Gamma)$. Let $U$ be the translation action of $\Gamma$ on $L^{2}(\nu)$.  We prove Lemma \ref{lem:TranslationExtremeErgodic} at the end of this section.

	\begin{lemma}\label{lem:BoundedTranslationErgodicity}
		 If $w\in L^2(\nu)$ is represented by an $L^{2}(\nu)$ Cauchy sequence of  \emph{uniformly} bounded functions $f_n\in \ell^\infty(\Gamma)$ and $w$ is $U$\text{-}invariant, then $w\sim_\nu c1_\Gamma$, where $c=\langle w, 1_\Gamma\rangle_\nu$.
	\end{lemma}
	
	\begin{proof}
		The hypothesis states that for some $\alpha>0$, $w$ belongs to the $L^2(\nu)$-equivalence class of a Cauchy sequence of functions $f_n:\Gamma\to \mathbb C$, where $|f_n(x)|\leq \alpha$ for all $x\in \Gamma$, and that $\lim_{n\to\infty}\|f_n-\gamma\cdot f_n\|_{\nu}=0$ for all $\gamma\in \Gamma$.  Splitting $f_n$ into real and imaginary parts $g_{n}$ and $h_n$, we see that $(g_n)_{n\in\mathbb{N}}$ and $(h_n)_{n\in\mathbb{N}}$ are themselves Cauchy, and represent $U$\text{-}invariant elements of $L^2(\nu)$.  It therefore suffices to prove the lemma under the assumption that the $f_n$ are real-valued.  By adding a constant and scaling, we may assume $f_n:\Gamma\to [0,1]$.
		
		Let $c=\langle w, 1_\Gamma \rangle_\nu$, meaning $c=\lim_{n\to\infty} \nu(f_n)$.  If $c=0$, then $\lim_{n\to\infty} \nu(|f_n|)=0$, meaning $w=0$ in $L^2(\nu)$.  Likewise, if $c=1$, it is easy to verify that $w\sim_\nu 1_\Gamma$.  So we assume $0<c<1$.  We define means $m_1$ and $m_2$ on $\ell^\infty(\Gamma)$ by $m_1(g) := \frac{1}{c}\langle g,w\rangle_\nu$, $m_2(g):=\frac{1}{1-c}\langle g, 1_\Gamma-w\rangle_\nu$.    Then $\nu=cm_1+(1-c)m_2$.   Since $w\in \mathcal H_{U\text{-inv}}$, we get that $m_{1}$ and $m_{2}$ are invariant means.  Since $\nu$ is extreme, we conclude that $m_1=m_2=\nu$, meaning $(1-c)\langle g,w\rangle_\nu = c\langle g,1_\Gamma-w\rangle_\nu$.  Simplifying, we get $\langle g,w\rangle_\nu = c\langle g,1_\Gamma\rangle_\nu=\langle g, c1_{\Gamma}\rangle_{\nu}$ for every $g\in \ell^\infty(\Gamma)$.  Since $\ell^{\infty}(\Gamma)$ forms a dense subspace of $L^{2}(\nu)$,  we see that $\langle v,w\rangle_\nu = \langle v, c1_\Gamma\rangle_\nu$ for all $v\in L^2(\nu)$.  This implies $w=c1_{\Gamma}$ in $L^{2}(\nu)$.
	\end{proof}

\subsection{Truncations}

Let $\alpha>0$, and define $\psi_{\alpha}:\mathbb C\to \mathbb C$ by
\[
\psi_{\alpha}(z) := \begin{cases}
	z &\text{ if } |z|<\alpha\\
	\alpha z/|z| &\text{ if } |z|\geq \alpha
\end{cases}
\]
Note that $\psi_{\alpha}$ is Lipschitz, with constant 1:  $|\psi_{\alpha}(z)-\psi_{\alpha}(z')|\leq |z-z'|$ for all $z, z'\in \mathbb C$.

For a function $f$ with codomain $\mathbb C$, define the truncation $f^{(\alpha)}$ by $f^{(\alpha)}:=\psi_{\alpha} \circ f$. 	Then $f^{(\alpha)}(x)=f(x)$ when $|f(x)|\leq \alpha$, and $|f^{(\alpha)}(x)|=\alpha$ when $|f(x)|\geq \alpha$.
	
\begin{lemma}\label{lem:Truncate}
	Let $m$ be a mean on $\ell^{\infty}(\Gamma)$,   let $w\in L^{2}(m)$, and let $(f_{n})_{n\in \mathbb N}$ be a Cauchy sequence of elements of $\ell^{\infty}(\Gamma)$ representing $w$.  
	
	For each $\alpha >0$, $(f_{n}^{(\alpha)})_{n\in \mathbb N}$ is a Cauchy sequence representing an element $w^{(\alpha)}\in L^{2}(m)$. 
	
	Furthermore $\lim_{\alpha\to\infty} w^{(\alpha)}=w$ in $L^{2}(m)$.
\end{lemma}
\begin{proof}
Let $\alpha>0$ and $f, g\in \ell^{\infty}(\Gamma)$.  Pointwise we have $|f^{(\alpha)}-g^{(\alpha)}|= |\psi_{\alpha}(f)-\psi_{\alpha}(g)|\leq |f-g|$, so 
\begin{equation}\label{eq:falpha}
	\|f^{(\alpha)} -g^{(\alpha)}\|_{m}\leq \|f-g\|_{m}.
\end{equation}  It follows that if $(f_{n})_{n\in \mathbb N}$ is Cauchy, so is $(f_{n}^{(\alpha)})_{n\in \mathbb N}$.  Thus $w^{(\alpha)}$ is well-defined.

Now let $(f_{n})_{n\in \mathbb N}$ be a Cauchy sequence of elements of $\ell^{\infty}(\Gamma)$ representing $w$.   We will prove that $\lim_{\alpha\to\infty}w^{(\alpha)}=w$.  Let $\varepsilon>0$.

Fix $n\in \mathbb N$ such that $\|f_{j}-f_{\ell}\|_{m}<\varepsilon$ for all $j,\ell \geq n$. Then for all $\alpha>0$ and all $j, \ell\geq n$, (\ref{eq:falpha}) implies $\|f_{j}^{(\alpha)}-f_{\ell}^{(\alpha)}\|_{m}<\varepsilon$. Choose $\alpha_{0}>0$ so that $f_{n}^{(\alpha_{0})}=f_{n}$. For every $\alpha\geq \alpha_{0}$ and all $j,\ell \geq n$, we have \[\|f_{j}^{(\alpha)}-f_{\ell}\|_{m}\leq \|f_{j}^{(\alpha)}-f_{n}^{(\alpha)}\|_{m}+\|f_{n}^{(\alpha)}-f_{n}\|_{m}+\|f_{n}-f_{\ell}\|_{m}<2\varepsilon,\]
so $\|w^{(\alpha)}-w\|_{m}<2\varepsilon$.  Since $\varepsilon>0$ was arbitrary, this shows that $\lim_{\alpha\to\infty} w^{(\alpha)}=w$. \end{proof}

\begin{proof}[Proof of Lemma \ref{lem:TranslationExtremeErgodic}]
To prove part	\ref{item:TranslationErgodicity}, assume $\nu\in \mathcal M_{\tau}^{ext}(\Gamma)$  and $w\in L^2(\nu)$ satisfies $U_\gamma w = w$ for all $\gamma \in \Gamma$.  Let $(f_n)_{n\in \mathbb N}$ be an $L^2(\nu)$-Cauchy sequence of elements of $\ell^\infty(\Gamma)$ representing $w$.    Let $\alpha>0$.  By Lemma \ref{lem:Truncate}, we get that $(f_{n}^{(\alpha)})$ is Cauchy, and therefore converges to an element $w^{(\alpha)}\in L^{2}(\nu)$. It is straightforward to verify that $w^{(\alpha)}$ is $U$-invariant, assuming $w$ is $U$-invariant.

  By Lemma \ref{lem:BoundedTranslationErgodicity}, $w^{(\alpha)}\sim_\nu c_\alpha 1_\Gamma$, where $c_\alpha =\langle w^{(\alpha)}, 1_\Gamma\rangle_\nu$.  Now Lemma \ref{lem:Truncate} implies $w=\lim_{\alpha\to\infty} c_{\alpha} 1_{\Gamma}=c 1_{\Gamma}$, where $c=\lim_{\alpha\to\infty} c_{\alpha}$.  This completes the proof of Part \ref{item:TranslationErgodicity}.

To prove \ref{item:TranslationEigenspace}, let $\chi\in\widehat{\Gamma}$, and suppose $w\in L^2(\nu)$ satisfies $U_\gamma w=\chi(\gamma)w$ for all $\gamma\in\widehat{\Gamma}$.  We will show that $w$ belongs to the subspace of $L^{2}(\nu)$ spanned by $\{\overline{\chi}\}$. Let $(f_n)_{n\in \mathbb N}\in \ell^{\infty}(\Gamma)$ be a Cauchy sequence representing $w$.  Then $(\chi f_{n})_{n\in \mathbb N}$ represents an element $v$ of $L^{2}(\nu)$.  We have $U_{\gamma}v=\lim_{n\to\infty} \overline{\chi(\gamma)}\chi U_{\gamma}f_{n}=\lim_{n\to\infty} \overline{\chi(\gamma)}\chi\chi(\gamma)f_{n}=\chi w=v$ for each $\gamma \in \Gamma$.  Thus $v$ is $U$\text{-}invariant, so part \ref{item:TranslationErgodicity} implies $v=c1_{\Gamma}$ for some $c\in \mathbb C$.  In other words, $\chi f_{n}$ converges to $c1_{\Gamma}$ in $L^{2}(\nu)$, so $f_{n}$ converges to $c\overline{\chi}$ in $L^{2}(\nu)$.  This shows that $w\sim_{\nu} c\overline{\chi}$, as desired.  \end{proof}
	\section{Proof of Lemma \ref{lem:Main}}\label{sec:ProofOfMain}

	Recall Definition \ref{def:meanFourierCoefficients}:  the Fourier coefficients of $f$ with respect to $\eta$ are $	\hat{f}^{\eta}(\chi):=\eta(f\overline{\chi}).$

\noindent Recall Lemma \ref{lem:Main}: 	let $\nu\in \mathcal M_{\tau}^{ext}(\Gamma)$, $\eta\in \mathcal M^{FS}(\Gamma)$, and $A$, $ B\subseteq \Gamma$.  Let $f=1_A$, $g=1_B$, and let $\tilde{f}^{(\nu)}$, $\tilde{g}^{(\eta)}$ be as in Lemma \ref{lem:DefBRN}.  Let $\tilde{h}=\tilde{f}^{(\nu)}*\tilde{g}^{(\eta)}$ and let $h=\tilde{h}|_{\Gamma}$.  Then the series
\begin{equation}\label{eq:RecallNuEtaConvolutionExpansion}
	\sum_{\chi\in\widehat{\Gamma}} \hat{f}^{\nu}(\chi)\hat{g}^{\eta}(\chi)\chi
\end{equation}
converges uniformly to $h$, and
\begin{enumerate}
	\item\label{item:RecallGeneralHAplusB} $\nu(h1_{A+B})=\nu(h)$;
	\item\label{item:RecallGeneralLevelContained} for all $\delta>0$, $\{h\geq \delta\}\subset_\nu A+B$;
	\item\label{item:RecallGeneralClopenSubset} if $V\subseteq \{\tilde{h}>0\}$ is compact, then $V\cap \Gamma \subset_{\nu} A+B$;
	\item\label{item:RecallGeneralMeasureOfAplusB} $\nu(A+B)\geq \mu_{b\Gamma}(\{\tilde{h}>0\})$.
\end{enumerate}

	\begin{proof}[Proof of Lemma \ref{lem:Main}]
Fix $A, B\subseteq \Gamma$, and let $f=1_{A}$, $g=1_{B}$.   Let $\tilde{f}=\tilde{f}^{(\nu)}$, $\tilde{g}=\tilde{g}^{(\eta)}$ be as in Lemma \ref{lem:DefBRN}.

Expanding $\tilde{h}:=\tilde{f}*\tilde{g}$ as a Fourier series, equation (\ref{eq:ConvolutionExpansion}) implies $\tilde{h}=\sum_{\chi\in\widehat{\Gamma}} \hat{\tilde{f}}(\tilde{\chi})\hat{\tilde{g}}(\tilde{\chi})\tilde{\chi}$ pointwise.    Lemma \ref{lem:DefBRN} allows us to replace the coefficients in the series and we get $\tilde{h}  =\sum_{\chi\in\widehat{\Gamma}} \hat{f}^{\nu}(\chi)\hat{g}^{\eta}(\chi)\tilde{\chi}$. This series converges uniformly: $\hat{f}^{\nu}(\chi)$ and $\hat{g}^{\eta}(\chi)$ belong to $\ell^{2}(\widehat{\Gamma})$ by Lemma \ref{lem:DefBRN}, so Cauchy-Schwarz implies $\sum_{\chi\in\widehat{\Gamma}} |\hat{f}^{\nu}(\chi)\hat{g}^{\eta}(\chi)|$ converges.   Hence $h:=\tilde{h}|_{\Gamma}=\sum_{\chi\in\widehat{\Gamma}} \hat{f}^{\nu}(\chi)\hat{g}^{\eta}(\chi)\chi$, and by Lemma \ref{lem:FStarGExpansion} we have $h\sim_{\nu} f\leftindex_{\nu}{*}_{\eta}\, g$.   This and Lemma \ref{lem:WeakSupportOnAplusB} imply Part \ref{item:GeneralHAplusB} of Lemma \ref{lem:Main}.  Part \ref{item:GeneralLevelContained} follows from Part \ref{item:GeneralHAplusB} and Lemma \ref{lem:EssentialLevelContainment}.
 
 To prove Part \ref{item:GeneralClopenSubset}, let $V\subseteq \{\tilde{h}>0\}$ be compact.  Let $\delta=(1/2)\min\{\tilde{h}(x):x\in V\}$.  Then $\delta>0$, since, $\tilde{h}$ is continuous, $V$ is compact, and $V\subseteq \{\tilde{h}>0\}$.  Thus $V\cap \Gamma\subseteq \{h\geq \delta\}$.  The essential containment $\{h\geq \delta\}\subset_{\nu} A+B$ from Part \ref{item:GeneralLevelContained} now implies $V\cap \Gamma \subset_{\nu} A+B$.

To prove Part \ref{item:GeneralMeasureOfAplusB}, let $\varepsilon>0$.  Inner regularity of Haar measure on compact groups provides a compact set $V\subseteq \{\tilde{h}>0\}$ with $\mu_{b\Gamma}(V)>\mu_{b\Gamma}(\{\tilde{h}>0\})-\varepsilon$.  Let $V'\subseteq \{\tilde{h}>0\}$ be a compact neighborhood of $V$, and let $\delta=\min\{\tilde{h}(x):x\in V'\}$, so that $\delta>0$. Let $\tilde{\phi}:b\Gamma\to [0,1]$ be a continuous function supported on $V'$ such that $\tilde{\phi}|_{V}=1$.  Let $\phi=\tilde{\phi}|_{\Gamma}$. Then 
\begin{align*}
	\nu(A+B) &\geq \nu(\{h\geq \delta\}) && \text{by part \ref{item:GeneralLevelContained}}\\ &\geq \nu(1_{V'\cap \Gamma}) && \text{since $V'\subseteq \{\tilde{h}>\delta\}$}\\ &\geq \nu(\phi) && \text{since } 1_{V'\cap \Gamma} \geq \phi \\ &=\int \tilde{\phi}\, d\mu_{b\Gamma} && \text{by Lemma \ref{lem:CorrectAveragingAndbGamma}} \\ &\geq \mu_{b\Gamma}(V) && \text{since } \tilde{\phi} \geq 1_{V}
	\\&\geq \mu_{b\Gamma}(\{\tilde{h}>0\})-\varepsilon.
\end{align*} 
Letting $\varepsilon\to 0$, we get $\nu(A+B)\geq \mu_{b\Gamma}(\{\tilde{h}>0\})$. \end{proof}

\section{Sumsets in level sets of convolutions, proof of Lemma \ref{lem:LevelSetSumset}}\label{sec:LevelSumset}
We will derive Lemma \ref{lem:LevelSetSumset} from the following special case.
\begin{lemma}\label{lem:PointsOfDensity}
	Let $K$ be a compact abelian group with Haar probability measure $\mu$, and let $A$, $B\subseteq K$ be $\mu$-measurable.  There are $\mu$-measurable sets $A' \subseteq A$, $B'\subseteq B$ such that $\mu(A')=\mu(A)$, $\mu(B')=\mu(B)$, and $1_A*1_B(x+y)>0$ for all $x\in A'$, $y\in B'$.  Consequently, $A'+B'\subseteq \{1_A*1_B>0\}$.
\end{lemma}

\begin{proof} Let $\mathcal U$ be a neighborhood base for $0_K$ consisting of compact symmetric neighborhoods of $0_K$, meaning $U=-U$ for every $U\in \mathcal U$. By Proposition 2.44 on of \cite[p.58]{Folland_CourseInAHA}, the collection $\{\psi_{U}:U\in \mathcal U\}$ of functions $\psi_U:=\mu(U)^{-1}1_{U}$ has the following property: for all $f\in L^2(\mu)$ and all $\varepsilon>0$, there is a neighborhood $V$ of $0_K$ such that for all $U\in \mathcal U$ with $U\subseteq V$, $\|\psi_{U}* f - f\|<\varepsilon$, where $\|\cdot\|$ denotes the $L^2(\mu)$ norm.

Let $f=1_A$ and $g=1_B$.	For each $n$, select $U_{n}\in \mathcal U$ such that
	\[\|\psi_{U_{n}}*f-f\|<1/n \qquad \text{ and } \qquad \|\psi_{U_{n}}*g-g\|<1/n.\]  Writing $\psi_{n}$  for $\psi_{U_{n}}$ and passing to a subsequence, we obtain 
	\begin{equation}\label{eq:ApproxIDae}
		\lim_{n\to\infty} \psi_{n}*f=f, \quad \lim_{n\to\infty} \psi_{n}*g=g \qquad \mu\text{-almost everywhere.}
	\end{equation}
	Let 
	\[A':=\{x\in K: \lim_{n\to\infty} \psi_{n}*f(x)=f(x)\}, \quad B':=\{y\in K:  \lim_{n\to\infty} \psi_{n}*g(y)=g(y)\}.\]	
	Let $x\in A'$ and $y\in B'$.  Note that $\psi_{n}*f(x)$ simplifies as $\mu(U_n)^{-1} \mu(A\cap(U_n+x))$, and similarly for $\psi_{n}*g(y)$.   So we may choose $n$ such that $\mu(A\cap (U_n+x))>0.6 \mu(U_n)$ and $\mu(B_n\cap (U_n+y))>0.6 \mu(U_n)$.  Since $U_{n}=-U_{n}$ and $\mu$ is invariant under the map $t\mapsto -t$, we have $\mu(B_{n}\cap (U_{n}+y))=\mu((B_{n}-y)\cap U_{n})=\mu((y-B_{n})\cap U_{n})$.  Similarly $\mu(A\cap (U_{n}+x))=\mu(U_{n}\cap (A-x))$. 	This means
	\begin{equation}\label{eq:90percent}
		\mu(U_n\cap (A-x)) > 0.6\mu(U_n),\quad \mu(U_n\cap (y-B)) > 0.6\mu(U_n).
	\end{equation}
	Then
	\begin{align*}
		1_A*1_B(x+y) &= \int 1_A(t)1_B(x+y-t)\, d\mu(t)\\
		&=  \int 1_A(t+x) 1_B(y-t)\, d\mu(t) && \\
		&= \int 1_{A-x}(t) 1_{B-y}(-t)\, d\mu(t)\\
		&= \int 1_{A-x}(t) 1_{y-B}(t)\, d\mu(t)\\
		&= \mu((A-x)\cap (y-B))\\
		&\geq \mu(U_n\cap (A-x)\cap (y-B))\\
		&\geq 0.2\mu(U_n),  && \text{by } (\ref{eq:90percent})
	\end{align*}
	so $1_A*1_B(x+y)>0$.
\end{proof}

Recall Lemma  \ref{lem:LevelSetSumset}:  	let $K$ be a compact abelian group with Haar probability measure $\mu$, and let $\tilde{f}$, $\tilde{g}:K\to [0,1]$ be $\mu$-measurable functions. Then there are $\mu$-measurable sets $\tilde{A}$, $\tilde{B}\subseteq K$ such that $\mu(\tilde{A})\geq \mu(\{\tilde{f}>0\})$, $\mu(\tilde{B})\geq \mu(\{\tilde{g}>0\})$, while  
$\tilde{f}(a)>0$, $\tilde{g}(b)>0$, and $\tilde{f}*\tilde{g}(a+b)>0$ for all $a\in \tilde{A}$, $b\in \tilde{B}$.

Consequently, $\{\tilde{f}*\tilde{g}>0\}$ is an open neighborhood of $\tilde{A}+\tilde{B}$.

\begin{proof}[Proof of Lemma \ref{lem:LevelSetSumset}]
	Let $\tilde{f}, \tilde{g}:K\to [0,1]$ be $\mu$-measurable functions with $\int\tilde{f}\, d\mu=\alpha$, $\int\tilde{g}\, d\mu=\beta$.
	
	For each $n$, let $C_n:=\{\tilde{f}>1/n\}$, $D_n:=\{\tilde{g}>1/n\}$.  By Lemma \ref{lem:PointsOfDensity} choose  $A_n\subseteq C_n$ and $B_n\subseteq D_n$ so that $\mu(A_n)=\mu(C_n)$, $\mu(B_n)=\mu(D_n)$, and $1_{C_n}*1_{D_n}(x+y)>0$ for all $x\in A_n$, $y\in B_n$.
	
	For $n\in \mathbb N$, let $A_n':= \bigcap_{k=n}^{\infty} A_k$,  and $B_n':= \bigcap_{k=n}^{\infty} B_k$ so that $A_n'\subseteq A_n$, $B_n'\subseteq B_n$ for each $n$, while $A_1'\subseteq A_2' \subseteq \cdots$ and $B_1'\subseteq B_2' \subseteq \cdots$.  Furthermore, $\mu(A_n')=\mu(A_n)$ and $\mu(B_n')=\mu(B_n)$ for each $n$.
	
	Let $\tilde{A}:=\bigcup_{n=1}^\infty  A_n'$, $\tilde{B}:=\bigcup_{n=1}^\infty B_n'$.  Then $\tilde{A}\subseteq \{\tilde{f}>0\}$, $\tilde{B}\subseteq \{\tilde{g}>0\}$, $\mu(\tilde{A})=\mu\{\tilde{f}>0\}$, and $\mu(\tilde{B})=\mu(\{\tilde{g}>0\})$.  
	
	To prove that $\tilde{A}+\tilde{B}\subseteq \{\tilde{f}*\tilde{g}>0\}$, let $a\in \tilde{A}$, $b\in \tilde{B}$.  Our definition of $C_{n}$ and $D_{n}$ implies $\tilde{f}*\tilde{g}(x)\geq n^{-2}1_{C_{n}}*1_{D_{n}}(x)$ for every $x\in K$ and all $n\in \mathbb N$, so it suffices to prove that $1_{C_{n}}*1_{D_{n}}(a+b)>0$ for some $n\in \mathbb N$. To see this, note that there is an $n\in \mathbb N$ such that $a\in A_{n}$ and $b\in B_{n}$, and this implies $1_{C_{n}}*1_{D_{n}}(a+b)>0$.
	
	Finally, note that $\tilde{f}*\tilde{g}$ is continuous, by Lemma \ref{lem:ConvolutionIsContinuous}, so $\{\tilde{f}*\tilde{g}>0\}$ is open.
\end{proof}

\section{Continuous measures on \texorpdfstring{$\widehat{\Gamma}$}{Gamma} and FS-means}\label{sec:Continuous}
In this section we prove Theorem \ref{th:WeakErgodic} and Lemma \ref{lem:RequiredProperties}.

\subsection{The weak ergodic theorem}

Recall Theorem \ref{th:WeakErgodic}: let $\Gamma$ be a discrete abelian group.	If $m\in \mathcal M_{\tau}(\Gamma)$, $U$ is a unitary action of $\Gamma$ on a Hilbert space $\mathcal H$, $v, w\in \mathcal H$, and $\phi_{v,w}(\gamma)=\langle v, U_{\gamma} w\rangle$ is the corresponding matrix coefficient, then $	m(\phi_{v,w})=\langle v, P_{U\text{-inv}} w\rangle.$
\begin{proof}[Proof of Theorem \ref{th:WeakErgodic}]
	Writing $w=w_{0}+w_{1}$, where $w_{1}=P_{U\text{-inv}}w$ and $w_{0}\perp \mathcal H_{U\text{-inv}}$, we have $\phi_{v,w}=\langle v, w_{1}\rangle+\phi_{v,w_{0}}$.  It suffices to prove that $m(\phi_{v,w_{0}})=0$.  Consider the linear functional $\Lambda:\mathcal H\to \mathbb C$ given by
	\[
	\Lambda(x) := \int \langle x, U_{\gamma }w_{0}\rangle \, dm(\gamma).
	\]
	By the Reisz representation theorem for Hilbert spaces, there is a unique $z\in \mathcal H$ such that $\Lambda(x)=\langle x,z\rangle$ for all $x\in \mathcal H$.  We will show that $z=0$, by showing that $z\in \mathcal H_{U\text{-inv}}$ and $z\in \mathcal H_{U\text{-inv}}^{\perp}$.
	
	To see that $z\in \mathcal{H}_{U\text{-inv}}$, fix $x\in \mathcal H$ and let $t \in \Gamma$.  Then \begin{align*}
		\langle x,U_{t} z\rangle &=\langle U_{t}^{-1} x, z\rangle && \text{since $U_{t}$ is unitary}\\
		&=\int \langle U_{t}^{-1} x, U_{\gamma} w_{0}\rangle\, dm(\gamma)\\
		&=\int \langle  x, U_{t+\gamma} w_{0}\rangle\, dm(\gamma)\\
		&= \int \langle x, U_{\gamma} w_{0}\rangle\, dm(\gamma)  && \text{since $m$ is invariant}\\
		&=\langle x, z\rangle.
	\end{align*}  Since $\langle x,z\rangle = \langle x, U_{t} z\rangle$ for all $x\in \mathcal H$ and all $t\in \Gamma$, we have $U_{t}z=z$ for all $t\in \Gamma$.
	
	To see that $z\in \mathcal{H}_{U\text{-inv}}^{\perp}$, let $x\in \mathcal{H}_{U\text{-inv}}$.  Then $\langle x, U_{\gamma} w_{0}\rangle = 0$ for every $\gamma\in \Gamma$, so	$\langle x,z\rangle=\int \langle x, U_{\gamma} w_{0}\rangle \, dm(\gamma)=0$. 	This shows that $z=0$, hence $\Lambda(x)=0$ for all $x\in \mathcal H$, meaning $m(\phi_{v,w_{0}})=0$.
\end{proof}

\subsection{Multiplication action}  
\begin{definition} Let $X$ be a topological space. A Borel measure $\mu$ on $X$ is \emph{outer regular} if $\mu(E)=\inf\{\mu(U):E\subset U, U \text{ is open}\}$ for every $\mu$-measurable $E\subseteq X$.  It is \emph{inner regular} if $\mu(E)=\sup \{\mu(V):V\subseteq E, V\text{ is compact}\}$.
	
	A \emph{Radon measure} on $X$ is a Borel measure that is finite on all compact subsets of $X$, outer regular on all Borel sets, and inner regular on open sets.
\end{definition}
In the setting where $\Gamma$ is a countable abelian group, $\widehat{\Gamma}$ is compact and metrizable, hence every Borel measure on $\widehat{\Gamma}$ is a Radon measure.  In general (when $\Gamma$ may be uncountable), we restrict our attention to Radon measures on $\widehat{\Gamma}$, as they retain the approximation properties we require.  See Chapter 7 of \cite{Folland_RealAnalysis} for details.

The next lemma is one of the building blocks spectral theory for unitary actions of discrete abelian groups.

\begin{lemma}\label{lem:Meigenspaces}  Let $\sigma$ be a positive finite Radon measure on $\widehat{\Gamma}$, and let $M$ be the action of $\Gamma$ on $L^{2}(\sigma)$ by unitary operators $M_{\gamma}$, where $(M_{\gamma} f)(\chi):=\chi(\gamma)f(\chi)$.  For each $\psi\in \widehat{\Gamma}$, the $\psi$-eigenspace $L^{2}(\sigma)_{\psi}$ of $M$ is at most one-dimensional, and is nontrivial if and only if $\sigma(\{\psi\})>0$.  If $\sigma(\{\psi\})>0$, then $L^{2}(\sigma)_{\psi}$ is spanned by $1_{\{\psi\}}$.
\end{lemma}

\begin{proof}
	It suffices to prove that if $\sigma(\{\psi\})>0$ then $1_{\{\psi\}}$ spans $L^{2}(\sigma)_{\psi}$.  So assume $f\in L^{2}(\sigma)_{\psi}$, meaning \begin{equation}\label{eq:chiGammaF}\chi(\gamma)f(\chi)=\psi f(\chi) \qquad \text{for } \sigma\text{-a.e. } \chi\in \widehat{\Gamma}, \text{ and all } \gamma\in \Gamma.
	\end{equation}  We will show that $f$ belongs to the span of $\{1_{\{\psi\}}\}$, by proving the following:
	\begin{equation}\label{eq:Target} \text{for every neighborhood } U \text{ of } \psi, f(\chi)=0 \text{ for } \sigma\text{-a.e.  }\chi\notin U.
	\end{equation}  Since $\sigma$ is a Radon measure, we may take a sequence of neighborhoods $U_{n}$ of $\psi$ with $\lim_{n\to\infty}\sigma(U_{n})=\sigma(\{\psi\})$ and conclude from (\ref{eq:Target}) that $f(\chi)=0$ for $\sigma$-almost every $\chi\neq \psi$.
	\begin{claim}
		For every  $\phi\in C(\widehat{\Gamma})$ and $\sigma$-almost every $\chi$, we have $\phi(\chi)f(\chi)=\phi(\psi)f(\chi)$.
	\end{claim}
	\begin{proof}[Proof of Claim]
		Let $\phi \in C(\widehat{\Gamma})$ and $\varepsilon>0$.  By Pontryagin duality, there are $\gamma_{1},\dots, \gamma_{k}\in \Gamma$ and $c_{1},\dots,c_{k}\in \mathbb C$ such that the trigonometric polynomial $p:\widehat{\Gamma}\to \mathbb C$ given by $p(\chi):=\sum_{j=1}^{k} c_{j}\chi(\gamma_{j})$ satisfies $\|\phi-p\|_{\infty}<\varepsilon$.  Then \[|\phi(\chi)f(\chi)-p(\chi)f(\chi)|\leq \varepsilon|f(\chi)|\] for all $\chi$, and (\ref{eq:chiGammaF}) implies  $p(\chi)f(\chi)=p(\psi) f(\chi)$ for $\sigma$-almost every $\chi$.  Equivalently,
		\begin{equation}\label{eq:Rearranged}
			|p(\chi)-p(\psi)||f(\chi)|=0  \qquad \text{for } \sigma\text{-a.e. } \chi. 
		\end{equation}Thus 
		\begin{align*}
			|&\phi(\chi)f(\chi)-\phi(\psi)f(\chi)|\\
			&=\big|\phi(\chi)-\phi(\psi)\big||f(\chi)|\\
			&\leq \big(|\phi(\chi)-p(\chi)|+|p(\chi)-p(\psi)|+|p(\psi)-\phi(\psi)|\big)|f(\chi)|\\
			&= \big|\phi(\chi)-p(\chi)\big||f(\chi)|+\big|p(\chi)-p(\psi)\big||f(\chi)|+\big|p(\psi)-\phi(\psi)\big||f(\chi)| \\
			&\leq \varepsilon |f(\chi)| + 0 + \varepsilon|f(\chi)| \qquad \text{ for } \sigma\text{-a.e. } \chi && \text{by (\ref{eq:Rearranged})}.
		\end{align*}
		Letting $\varepsilon\to 0$, we get $\phi(\chi)f(\chi)=\phi(\psi)f(\chi)$ for $\sigma$-almost every $\chi$. \end{proof}
	To prove (\ref{eq:Target}), let $U$ be a neighborhood of $\psi$, and let $\phi:\widehat{\Gamma}\to [0,1]$ be a continuous function with $\phi(\psi)=1$ which vanishes outside $U$.  The claim then implies $\phi(\chi)f(\chi)=\phi(\psi)f(\chi)=f(\chi)$ for $\sigma$-almost every $\chi$, which means $f(\chi)=0$ for $\sigma$-a.e.~$\chi\notin U$.
\end{proof}

\subsection{Fourier transforms of measures} For a finite positive Radon measure $\sigma$ on $\widehat{\Gamma}$, its \emph{Fourier transform} is $\hat{\sigma}:\Gamma\to \mathbb C$, given by $\hat{\sigma}(\chi):=\int \chi(\gamma)\, d\sigma(\chi)$.  The following is the Bochner-Herglotz theorem for this setting,  with some additional detail regarding eigenspaces.
\begin{theorem}\label{th:BochnerHerglotz}
	Let $U$ be a unitary action of $\Gamma$ on a Hilbert space $\mathcal H$, let $w\in \mathcal H$, and let $\phi_{w,w}:\Gamma \to \mathbb C$,  $\phi_{w,w}(\gamma):=\langle w,U_{\gamma} w\rangle$ be the corresponding matrix coefficient.  Then there is a finite Radon measure on $\widehat{\Gamma}$ such that
	\begin{equation}\label{eq:Bochner}	
		\phi_{w,w}(\gamma) = \hat{\sigma}(\gamma) \qquad \text{ for all } \gamma \in \Gamma.
	\end{equation}
	Furthermore, for each $\chi\in \widehat{\Gamma}$, the $\chi$-eigenspace of $U$ is trivial if and only if $\sigma(\{\chi\})=0$.
\end{theorem}
\begin{proof}
	The existence of the Radon measure $\sigma$ satisfying (\ref{eq:Bochner}) is given by Theorem 4.19 in \cite[p.~103]{Folland_CourseInAHA}. 
	
	To prove the second assertion, consider the unitary action $M$ of $\Gamma$ on $L^{2}(\sigma)$, given by $(M_{\gamma}f)(\chi)=\chi(\gamma)f(\chi)$. We exhibit the (well-known) unitary equivalence $\Phi$ between $M$ and the restriction of $U$ to $\mathcal H^{(w)}$ ($=$ the closed subspace spanned by $\{U_{\gamma}w:\gamma \in \Gamma\}$).  For finite linear combinations $\sum c_{\gamma}U_{\gamma}w$, define $\Phi(\sum c_{\gamma}U_{\gamma}w):=\sum c_{\gamma}M_{\gamma}1_{\widehat{\Gamma}}$.  It is easy to verify that $\Phi$ so defined is a linear isometry from a dense subspace of $\mathcal H^{(w)}$ into $L^{2}(\sigma)$, so $\Phi$ extends to a linear isometry defined on all of $\mathcal H^{(w)}$. By Pontryagin duality, every continuous function $f$ on $\widehat{\Gamma}$ can be uniformly approximated by a linear combination of functions of the form $\chi\mapsto \chi(\gamma)$.  Since $\sigma$ is a Radon measure, the continuous functions on $\widehat{\Gamma}$ form a dense subspace of $L^{2}(\sigma)$ (Proposition 7.9 in Ch.7 of \cite{Folland_RealAnalysis}).  Thus the image of $\Phi$ is dense in $L^{2}(\sigma)$.  Since $\Phi$ is an isometry, this implies $\Phi$ is onto.  	Thus $\Phi$ is a unitary equivalence between the unitary actions $U$ and $M$.  Fixing $\chi\in \widehat{\Gamma}$, Lemma \ref{lem:Meigenspaces} says that the $\chi$-eigenspace of $M$ is spanned by $1_{\{\chi\}}$, the characteristic function of the singleton $\{\chi\}$.  Thus the $\chi$-eigenspace of $U$ restricted to $\mathcal H^{(w)}$ is nontrivial if and only if $\sigma(\{\chi\})>0$.
\end{proof}

\subsection{Polarization identity}

Let $U$ be a unitary action of $\Gamma$ on a Hilbert space $\mathcal H$, and let $v,w \in \mathcal H$. Setting $z_1= v+w$, $z_2=v-w$, $z_3=v+iw$, and $z_4=v-iw$, one may verify
\begin{equation}\label{eq:polarization}
	4\langle v,  U_\gamma w\rangle = \langle  z_1, U_\gamma z_1\rangle - \langle z_2,U_\gamma  z_2\rangle + i\langle  z_3,U_\gamma  z_3\rangle - i\langle  z_4,U_\gamma z_4\rangle.
\end{equation}

For the next lemma, recall that if $U$ is a unitary action on a Hilbert space $\mathcal H$, then $\mathcal H_{c}$ is the closure of the span of the $U$-eigenspaces $\mathcal H_{\chi}$; thus $w\in \mathcal H_{c}^{\perp}$ is equivalent to $w$ being orthognal to every eigenspace of $U$. 

\begin{lemma}\label{lem:AnnihilatesImpliesFS}
	Let $\eta$ be a mean on $\ell^{\infty}(\Gamma)$ such that $\eta(|\hat{\sigma}|^{2})=0$ for every continuous positive finite Radon measure on $\widehat{\Gamma}$.  If $U$ is a unitary action of $\Gamma$ on a Hilbert space $\mathcal H$ and $w\in \mathcal H_{c}^{\perp}$, then $\eta(|\phi_{v,w}|^{2})=0$ for all $v\in \mathcal H$.	
\end{lemma}

\begin{proof}
	Let $U$ be a unitary action of $\Gamma$ on $\mathcal H$, and assume $w\in \mathcal H_{c}^{\perp}$.  
	
	Writing $v=v_{0}+v_{1}$, where $v_{1}\in \mathcal H_{c}$, we have $\phi_{v,w}=\phi_{v_{1},w}+\phi_{v_{0},w}$.  
	Since $\mathcal H_{c}$ and $\mathcal H_{c}^{\perp}$ are mutually orthogonal, closed, and invariant under $U$, we get that $\phi_{v_{1},w}=0$. It therefore suffices to prove the statement under the additional assumption that $v\in \mathcal H_{c}^{\perp}$.
	
	The polarization identity allows us to write $|\phi_{v,w}|^{2}$ as a linear combination of functions the form $\gamma \mapsto \langle z_{1}, U_\gamma z_{1}\rangle\overline{\langle z_{2}, U_\gamma z_{2}\rangle}$, where each $z_{i}\in \mathcal H_{c}^{\perp}$.  Since $\mathcal H_{c}^{\perp}$ has no nontrivial eigenspaces, Theorem \ref{th:BochnerHerglotz} provides a continuous Radon measure $\sigma_{i}$ on $\widehat{\Gamma}$  such that $\langle z_{i}, U_{\gamma} z_{i}\rangle=\hat{\sigma}_{i}(\gamma)$ for each $\gamma\in \Gamma$.  By hypothesis, we have $\eta(|\hat{\sigma}_{i}|^{2})=0$ for $i=1,2$. Applying Cauchy-Schwarz, we have $|\eta(\hat{\sigma}_{1}\overline{\hat{\sigma}_{2}})|\leq \eta(|\hat{\sigma}_{1}|^{2})\eta(|\hat{\sigma}_{2}|^{2})=0$.
\end{proof}

\subsection{Proof of Lemma \ref{lem:RequiredProperties}}

Lemma \ref{lem:RequiredProperties} is  \ref{item:CorrectAverage}$\iff$\ref{item:BothHilbert} in the following lemma.
\begin{lemma}\label{lem:FSequivalents}
	Let $\eta$ be a mean on $\ell^{\infty}(\Gamma)$.  The following are equivalent.
	
	\begin{enumerate}
		\item\label{item:CorrectAverage} For every unitary action $U$ of $\Gamma$ on a Hilbert space $\mathcal H$ and all $v, w\in \mathcal H$, we have $\eta(\phi_{v,w})=\langle v, P_{U\text{-inv}} w\rangle$.
		\item\label{item:AtomAt0} For every  positive finite Radon measure $\sigma$ on $\widehat{\Gamma}$, we have $\eta(\hat{\sigma})=\sigma(\{\chi_{0}\})$, where $\chi_{0}\in\widehat{\Gamma}$ is the trivial character.
		
		\item\label{item:StandardFS} Both of the following hold: 
		\begin{enumerate}
			\item for every continuous positive finite Radon measure $\sigma$ on $\widehat{\Gamma}$, $\eta(|\hat{\sigma}|^{2})=0$
			
			\item  $\eta(\chi)=0$ for every nontrivial $\chi\in\widehat{\Gamma}$.
		\end{enumerate}
		\item\label{item:BothHilbert} Both of the following hold: 
		
		\begin{enumerate}
			\item for every action $U$ of $\Gamma$ by unitary operators on a Hilbert space $\mathcal H$, and all $w\in \mathcal H_{c}^{\perp}$,  $\eta(|\phi_{v,w}|^{2})=0$ for all $v\in \mathcal H$.
			\item\label{eq:BothHilbertDiscrete}  $\eta(\chi)=0$ for every nontrivial $\chi \in \widehat{\Gamma}$.
		\end{enumerate}
	\end{enumerate}
\end{lemma}

\begin{proof}
	\ref{item:CorrectAverage}$\implies$\ref{item:AtomAt0}.	Assume \ref{item:CorrectAverage}, and let $\sigma$ be a positive finite Radon measure on $\widehat{\Gamma}$.  Let $v=w=1_{\widehat{\Gamma}}\in L^{2}(\sigma)$.  Let $U$ be the multiplication action on $L^{2}(\sigma)$ given by  $(U_{\gamma }f)(\chi)=\chi(\gamma)f(\chi)$.  By Lemma \ref{lem:Meigenspaces}, we have $P_{U\text{-inv}}w=1_{\{\chi_{0}\}}$, where $\chi_{0}\in \widehat{\Gamma}$ is the trivial character.  Under assumption \ref{item:CorrectAverage} we have $\eta(\phi_{w,w})= \langle 1_{\widehat{\Gamma}}, 1_{\{\chi_{0}\}}\rangle_{L^{2}(\sigma)} = \sigma(\{\chi_{0}\})$.

	\ref{item:AtomAt0}$\implies$\ref{item:StandardFS}.	Let $\sigma$ be a continuous positive finite Radon measure on $\widehat{\Gamma}$.  Note that $|\hat{\sigma}(\gamma)|^{2}=\widehat{\sigma*\sigma}(\gamma)$, where $\sigma*\sigma$ is the Radon measure on $\widehat{\Gamma}$ given by
	\[
	\int f\, d\sigma*\sigma = \int \int f(\chi\overline{\psi})\, d\sigma(\chi)\, d\sigma(\psi).\] To verify that  $\sigma*\sigma(\{\chi_{0}\})=0$, evaluate $\int 1_{\{\chi_{0}\}}\, d \sigma*\sigma= \int \int 1_{\{\chi_{0}\}}(\chi \overline{\psi})\, d\sigma(\chi)\, d\sigma(\chi)$, and note that the inner integral is always $0$ since $\sigma$ is continuous.  To see that $\sigma*\sigma$ is a Radon measure, apply the Reisz representation theorem to the positive linear functional $f\mapsto \int \int f(\chi\overline{\psi})\, \sigma(\chi)\, d\sigma(\psi)$ on $C(\widehat{\Gamma})$.  The resulting measure agrees with $\sigma*\sigma$ on Borel sets.
	
	Now $\eta(|\hat{\sigma}|^{2})=\eta(\widehat{\sigma*\sigma})=\sigma*\sigma (\{\chi_{0}\})=0$.
	
	To see that $\eta(\chi)=0$ for every nontrivial $\chi\in\widehat{\Gamma}$, note that $\chi=\hat{\sigma}$, where $\sigma=\delta_{\chi}$, the Dirac mass at $\chi$.  When $\chi$ is nontrivial, then $\eta(\chi)=\sigma(\{\chi_{0}\})=\delta_{\chi}(\{\chi_{0}\})=0$.
	
	\ref{item:StandardFS}$\implies$\ref{item:BothHilbert} is immediate from Lemma \ref{lem:AnnihilatesImpliesFS}.

	To prove \ref{item:BothHilbert}$\implies$\ref{item:CorrectAverage},	assume \ref{item:BothHilbert}, let $U$ be a unitary action of $\Gamma$ on a Hilbert space $\mathcal H$, and let $v, w\in \mathcal H$.  Note that $\mathcal H_{U\text{-inv}} = \mathcal H_{\chi_{0}}\subseteq \mathcal H_{c}$.  Writing $w=w_{0}+w_{1}$ where $w_{0}\in \mathcal H_{c}^{\perp}$ and $w_{1}\in \mathcal H_{c}$, we then have $P_{U\text{-inv}}w_{1}=P_{U\text{-inv}} w$.
	
	Now $\phi_{v,w}=\phi_{v,w_{0}}+\phi_{v,w_{1}}$, and assumption \ref{item:BothHilbert} implies $\eta(\phi_{v,w_{0}})=0$. To verify that \begin{equation}\label{eq:etaphi} \eta(\phi_{v,w_{1}})=\langle v, P_{U\text{-inv}}w_{1}\rangle
	\end{equation} for all $v\in \mathcal H$,  fix such a $v$.  Expanding $w_{1}$ as a combination of mutually orthogonal eigenvectors $\sum_{\chi \in \widehat{\Gamma}} c_{\chi}w_{\chi}$, where $U_{\gamma} w_{\chi}=\chi(\gamma) w_{\chi}$ for all $\gamma\in \Gamma$, we have $P_{U\text{-inv}} w_{1}= c_{\chi_{0}}w_{\chi_{0}}$.  It therefore suffices to prove that $\int \langle v,U_{\gamma} w_{\chi}\rangle d\eta(\gamma)=0$ for all nontrivial $\chi$.  This follows from assumption (\ref{eq:BothHilbertDiscrete}) and the simplification $\int \langle v, U_{\gamma} w_{\chi}\rangle\, d\eta(\gamma) = \int \overline{\chi}(\gamma)\langle v,w_{\chi}\rangle \, d\eta(\gamma) =  \langle v,w_{\chi}\rangle \int \overline{\chi}(\gamma) \, d\eta(\gamma)=\eta(\overline{\chi})=0$.

	We now have $\eta(\phi_{v,w})=\eta(\phi_{v,w_{1}})=\langle v,P_{U\text{-inv}} w\rangle$, as desired. \end{proof}

\section{Hartman uniform distribution and FS-means}\label{sec:HartmanUD}
As mentioned in \S\ref{sec:FSmeans}, every invariant mean on $\ell^{\infty}(\Gamma)$ is an FS-mean.  Since abelian groups are amenable, we have examples of FS-means in every infinite abelian group.  In countable abelian groups there are many interesting examples of non\text{-}invariant FS-means.  However, the author is presently unaware of interesting examples of non\text{-}invariant FS-means in \emph{uncountable} abelian groups.

We discuss examples of FS-means in countable abelian groups and ask whether similar examples can be obtained in uncountable groups.

\begin{definition}\label{def:HartmanUD}
	Let $\Gamma$ be a discrete abelian group.  A sequence $(\gamma_{j})_{j\in \mathbb N}$ of elements of $\Gamma$ is \emph{Hartman uniformly distributed} (``Hartman-u.d.'') if for every nontrivial character $\chi\in\widehat{\Gamma}$, we have $\lim_{n\to\infty}\frac{1}{n}\sum_{j=1}^{n}\chi(\gamma_{j})=0$.  More generally, a sequence of \emph{countably additive} measures $\mu_{n}$ on $\mathcal P(\Gamma)$ is Hartman-u.d.~if $\lim_{n\to\infty} \int \chi \, d\mu_{n}=0$ for every nontrivial character $\chi\in\widehat{\Gamma}$.  
\end{definition}
Thus $(\gamma_{j})_{j\in \mathbb N}$ is Hartman-u.d. if and only if the sequence of measures $\mu_{n}:=\frac{1}{n}\sum_{j=1}^{n}\delta_{\gamma_{j}}$ is Hartman-u.d.

Specializing to $\mathbb Z$, we see that a sequence of integers $a_{j}$ is Hartman-u.d.~if and only if $\lim_{n\to\infty} \frac{1}{n}\sum_{j=1}^{n} e^{i a_{j}t}=0$ for all $t\in (0,2\pi)$.  Examples of Hartman-u.d.~sequences include the following. See Section 3 of \cite{BKQW} or \cite{Nair_PoincareI} for details and more examples.
\begin{enumerate}
	\item[$\bullet$] $a_{n}=\lfloor n^{5/2}\rfloor$.
	
	\item[$\bullet$] Fix $c>0$, $c\notin \mathbb N$.  Then $b_{n}=\lfloor n^{c}\rfloor$ is Hartman-u.d.
	
	\item[$\bullet$] Let $p(x)=\alpha_{k}x^{k}+\alpha_{k-1}x^{k-1}+\cdots + \alpha_{1}x_{1}+\alpha_{0}$ be a polynomial with real coefficients. If $\alpha_{r}/\alpha_{s}$ is irrational for some $1\leq r<s\leq k$, Theorem 1 of \cite{Veech_Well} says that $c_{n}=\lfloor p(n)\rfloor$ is Hartman-u.d.
	
	\item[$\bullet$] An example growing faster than any polynomial is discussed in \S1 of \cite{Boshernitzan_Homogeneous}. Let $1<c<3/2$.  Then $d_{n}=\lfloor e^{(\log n)^{c}}\rfloor$ is Hartman-u.d.~The values $c$ greater than or equal to $3/2$ for which such sequences are Hartman-u.d.~are unknown.
\end{enumerate}

\begin{remark}
	The term ``Hartman uniformly distributed'' is common, but not quite standard.  For example, \cite{Veech_Well} uses ``uniformly distributed,'' \cite{Boshernitzan_Homogeneous} uses ``homogeneously distributed,'' \cite{Griesmer_SumsetsDenseSparse} uses ``ergodic averaging sequence,'' and \cite{BergelsonFerreMoragues_ErgodicCorrespondence} uses ``ergodic sequence.''
\end{remark}

The following well-known lemma will produce FS-means from Hartman-u.d. sequences.  

\begin{lemma}\label{lem:HartmanUDisFS}
	If $\mu_{j}$ is a Hartman-u.d.~sequence of measures on $\mathcal P(\Gamma)$ and $\sigma$ is a positive finite Radon measure on $\widehat{\Gamma}$, then $\lim_{n\to\infty}\int \hat{\sigma}(\gamma)\, d\mu_{n}(\gamma)=\sigma(\{\chi_{0}\})$.
\end{lemma}

\begin{proof}
	With $\mu_{j}$ and $\sigma$ as in the hypothesis, we have
	\begin{equation}\label{eq:Fub}
		\int \hat{\sigma}(\gamma)\, d\mu_{n}(\gamma)	= \int \int  \chi(\gamma) \, d\sigma(\chi)  \, d\mu_{n}(\gamma)=\int \int \chi(\gamma)\, d\mu_{n}(\gamma) \, d\sigma(\chi), 
	\end{equation}
	where we applied Fubini's theorem to get the second equality -- this is where we use the countable additivity of $\mu_{n}$.  Taking the limit inside the integral in the right-hand side of (\ref{eq:Fub}), we get $\int \chi(\gamma)\, d\mu_{n}(\gamma)=0$ for all $\chi$ except $\chi_{0}$.  The dominated convergence theorem then implies that the limit in (\ref{eq:Fub}) simplifies to $\int 1_{\{\chi_{0}\}}\, d\sigma(\chi)$, which is nothing but $\sigma(\{\chi_{0}\})$, as desired.   
\end{proof}

\begin{remark}
	Note that in applying the dominated convergence theorem, we require that $(\mu_{n})_{n\in \mathbb N}$ be a sequence, rather than a net.  This leads to Question \ref{q:HartmanMean}.
\end{remark}

Thus every Hartman-u.d.~sequence of elements of $\Gamma$ produces an FS-mean $\eta$ on $\ell^{\infty}(\Gamma)$: if $(a_{n})_{n\in \mathbb N}$ is Hartman-u.d., we apply Lemma \ref{lem:MeanFromSequence} to get a mean $\eta$ satisfying $\eta(f)=\lim_{N\to\infty} \frac{1}{N}\sum_{n=1}^{N} f(a_{n})$ for all $f\in \ell^{\infty}(\Gamma)$ where the limit exists.  In particular, Lemma \ref{lem:HartmanUDisFS} implies $\eta(\hat{\sigma})=\sigma(\{\chi_{0}\})$ for every positive finite Radon measure on $\widehat{\Gamma}$.  By Lemma \ref{lem:FSequivalents}, we get that $\eta$ is an FS-mean.

Since $a_{n}=\lfloor n^{5/2}\rfloor$ is a Hartman-u.d.~sequence, this shows that there is is an FS-mean $\eta$ on $\mathbb Z$ such that $\eta(\{\lfloor n^{5/2} \rfloor:n\in \mathbb N\})=1$.

\begin{question}\label{q:HartmanMean}
	Call a mean $\eta$ on $\ell^{\infty}(\Gamma)$ \emph{Hartman-u.d.}~if $\eta(\chi)=0$ for every nontrivial $\chi\in\widehat{\Gamma}$. Is every Hartman-u.d.~mean on $\ell^{\infty}(\Gamma)$ also an FS-mean?
\end{question}

\section{Finite index subgroups}

\begin{lemma}\label{lem:FiniteIndex}
		Let $\Gamma$ be a discrete abelian group and $b\Gamma$ its Bohr compactification.  Let $\tilde{K}\leq  b\Gamma$ is a $\mu_{b\Gamma}$-measurable finite index subgroup and let $K=\tilde{K}\cap \Gamma$. Then
		
			\begin{enumerate}
			\item[(i)]  $K$ has finite index in $\Gamma$, and $[\Gamma:K]=[b\Gamma:\tilde{K}]$.
			
			\item[(ii)] If $x+\tilde{K}$ is a coset of $\tilde{K}$, then $(x+\tilde{K})\cap \Gamma$ is a coset of $K$.
			
			\item[(iii)] If $C\subseteq \Gamma$, then $C+K = (C+\tilde{K})\cap \Gamma$.
			
			\item[(iv)] The map $\rho: \Gamma/K\to b\Gamma/\tilde{K}$ given by $\rho(\gamma+K)=\gamma+\tilde{K}$ is a group isomorphism; its inverse is given by $\rho^{-1}(x+\tilde{K})=(x+\tilde{K})\cap \Gamma$.
			
			\item[(v)] If $C\subseteq b\Gamma$ and $\eta$ is an FS-mean on $\ell^{\infty}(\Gamma)$, then $\eta((C+\tilde{K})\cap \Gamma)=\mu_{b\Gamma}(C+\tilde{K})$.
		\end{enumerate}
	\end{lemma}

	\begin{proof}
(i) Assuming $\tilde{K}$ is $\mu_{b\Gamma}$-measurable and has finite index in $b\Gamma$, we get that $\mu_{b\Gamma}(\tilde{K})>0$.  Hence Theorem \ref{th:SteinhausCompactAbelian} implies $\tilde{K}-\tilde{K}$ has nonempty interior.  Since $\tilde{K}$ is a subgroup, we have $\tilde{K}=\tilde{K}-\tilde{K}$, so $\tilde{K}$ has nonempty interior. Thus $\tilde{K}-\tilde{K}$ contains a neighborhood of the identity.  Then $\tilde{K}$ contains a neighborhood of the identity, so $\tilde{K}$ is open.

Let $d$ be the index of $\tilde{K}$ in $b\Gamma$. Fix coset representatives $x_{1},\dots,x_{d}$ of $\tilde{K}$, and choose a neighborhood $V$ of the identity in $b\Gamma$ such that $x_{j}+V\subseteq x_{j}+\tilde{K}$ for each $j$.  Since $\Gamma$ is dense in $b\Gamma$, we may choose $\gamma_{j}\in (x_{j}+V)\cap \Gamma$ for each $j$.  Then $\gamma_{j}+\tilde{K}=x_{j}+\tilde{K}$ for each $j$.  We claim that $\Gamma = \bigcup_{j\leq d} \gamma_{j} + K.$  Fix $\gamma\in \Gamma$ and $j\leq d$ such that $\gamma \in \gamma_{j}+\tilde{K}$.  Then $\gamma-\gamma_{j}\in \tilde{K}$ and $\gamma-\gamma_{j}\in \Gamma$, so $\gamma-\gamma_{j}\in K$.  Then $\gamma\in \gamma_{j}+K$.  Thus $\Gamma$ is covered by at most $d$ cosets of $K$.  The disjointness of the collection $\{\gamma_{j}+\tilde{K}:j\leq d\}$ implies that the $\gamma_{j}+K$ are mutually disjoint, as well, so $[\Gamma:K]=d$.

(ii) Let $x\in b\Gamma$.  As in part (i), we may find a $\gamma\in \Gamma$ such that $x+\tilde{K}=\gamma+\tilde{K}$.  Then $(x+\tilde{K})\cap \Gamma = (\gamma+ \tilde{K})\cap \Gamma$, and we see that $y\in (\gamma+\tilde{K})\cap \Gamma \iff y-\gamma \in \tilde{K}\cap \Gamma=K$.  So $y\in \gamma+K \iff y\in (\gamma +\tilde{K})\cap \Gamma$.

(iii) It suffices to prove that if $\gamma\in \Gamma$, then $\gamma+K=(\gamma+\tilde{K})\cap \Gamma$.  The containment $\gamma+K\subseteq (\gamma+\tilde{K})\cap \Gamma$ is clear. To prove the reverse containment, let $y\in (\gamma+\tilde{K})\cap \Gamma$.  Then $y=\gamma + \tilde{k} \in \Gamma$, for some $\tilde{k}\in \tilde{K}$.  So $\tilde{k} = y-\gamma \in \tilde{K}$, meaning $y-\gamma \in K$.  Thus $y\in \gamma + K$.

(iv) To see that $\rho(\gamma+K):=\gamma+\tilde{K}$ defines an isomorphism, first note that it defines a homomorphism, since $\rho(\gamma+K)=\gamma+K+\tilde{K}$, and $\tilde{K}$ is a subgroup of $b\Gamma$ containing $K$.  Furthermore $\rho$ is surjective, since $\Gamma$ is dense in $b\Gamma$ and $\tilde{K}$ is open.  We see that $\rho$ is an isomorphism, being a surjective homomorphism between finite groups of the same cardinality.  To see that $\rho^{-1}(x+\Gamma)=(x+\tilde{K})\cap \Gamma$, fix $\gamma\in \Gamma$ so that $x+\tilde{K}=\gamma+\tilde{K}$ and $(x+\tilde{K})\cap \Gamma= \gamma +K$ (as part (ii) allows). Then $\rho(\gamma+K)= \gamma+\tilde{K}= x+\tilde{K}$.

(v)  As seein in part (i), $\tilde{K}$ is open, and therefore clopen. Then $1_{\tilde{K}}$ is continuous, so $\eta(K)=\eta(1_{\tilde{K}}|_{\Gamma})=\int 1_{\tilde{K}}\, d\mu_{b\Gamma}=\mu_{b\Gamma}(\tilde{K})$, by Lemma \ref{lem:CorrectAveragingAndbGamma}. Likewise $\eta(\gamma+K)=\mu_{b\Gamma}(\tilde{K})$ for every $\gamma\in \Gamma$.

We then have $\eta(C+K)=r\cdot \mu_{b\Gamma}(\tilde{K})$, where $r$ is the number of cosets of $K$ occupied by $C+K$. So it suffices to verify that the number of cosets of $K$ occupied by $(C+\tilde{K})\cap \Gamma$ equals the number of cosets of $\tilde{K}$ occupied by $C+\tilde{K}$.  Writing $C+\tilde{K}$ as a disjoint union $(\gamma_{1}+\tilde{K})\cup \dots \cup (\gamma_{r}+\tilde{K})$, part (iv) implies $\rho^{-1}(C+\tilde{K})=(\gamma_{1}+K)\cup \dots \cup (\gamma_{r}+K)$, as desired. \end{proof}

	\section{Large subsets of abelian groups}\label{sec:LargeSubsets}
Fix a discrete abelian group $\Gamma$ for this section.
	
	\subsection{Thickness and syndeticity}
	Let $A\subseteq \Gamma$.  We say that $A$ is 
\begin{enumerate}
		\item[$\bullet$] \emph{thick} if for all finite sets $F\subseteq \Gamma$, there is a $t\in \Gamma$ such that $F+t\subseteq A$.
		
		\item[$\bullet$] \emph{syndetic} if there is a finite set $F\subseteq \Gamma$ such that $F+A=\Gamma$.  
		
		\item[$\bullet$] \emph{piecewise syndetic} if there is a finite set $F\subseteq \Gamma$ such that $F+A$ is thick.
	\end{enumerate}

	\subsection{
		Finite embeddability} Given $A,B\subseteq \Gamma$, we write $A\prec B$ if for all finite $F\subseteq A$, there is a $t\in \Gamma$ such that $F+t\subseteq B$.

	\begin{lemma}\label{lem:UBDequivalents}
		Let $A\subseteq \Gamma$ and $\alpha>0$.  The following are equivalent.
		
		\begin{enumerate}
			\item\label{item:UBD} There is an invariant mean $m$ on $\ell^\infty(\Gamma)$ such that $m(1_A)\geq \alpha$.
			\item\label{item:FiniteDensity} For all finite $F\subseteq \Gamma$, there is a $t\in \Gamma$ such that $|(F+t)\cap A|\geq \alpha|F|$.
			\item\label{item:FolnerDensity} There is a F{\o}lner net $\mb \Phi$ such that $d_{\mb \Phi}(A)\geq \alpha$.
		\end{enumerate}
		Assuming $\Gamma$ is countable, \ref{item:FolnerDensity} is equivalent to the condition that $d_{\mb\Phi}(A)\geq \alpha$ for some F{\o}lner sequence $\mb \Phi$.
	\end{lemma}
	
	\begin{proof} \ref{item:UBD} $\implies$ \ref{item:FiniteDensity}: let $m$ be an invariant mean with $m(1_A)\geq \alpha$, and let $F\subseteq \Gamma$ be finite.  Define $\phi:\Gamma\to \mathbb R$ by $\phi(t):= |(F+t)\cap A|$.  Then $\phi$ can be simplified as $\sum_{\gamma\in \Gamma} 1_{F+t}(\gamma)1_A(\gamma)=\sum_{\gamma \in F} 1_{A-t}(\gamma)=\sum_{\gamma\in F} 1_{A-\gamma}(t)$.  Since $m$ is invariant, we have $m(1_{A-\gamma})=m(1_A)$ for every $\gamma$.  Thus $m(\phi)=m(1_A)|F|\geq \alpha |F|$.  Then for all $\varepsilon>0$, there is a $t\in \Gamma$ such that $\phi(t)\geq \alpha |F|-\varepsilon$.  Since $\phi$ is integer valued, there must be some $t$ with $\phi(t)\geq |F|$.
		
		\ref{item:FiniteDensity} $\implies$ \ref{item:FolnerDensity}.  Assume \ref{item:FiniteDensity} holds, and let $(\Phi_j)_{j\in I}$ be any F{\o}lner net.  For each $j\in I$, apply \ref{item:FiniteDensity} to choose $t_j$ so that $|(\Phi_j+t_j)\cap A|\geq \alpha |\Phi_j|$.  Then $\mb \Phi'=(\Phi_j+t_j)_{j\in I}$ is a F{\o}lner net with $d_{\Phi'}(A)\geq \alpha$.
		
		To prove \ref{item:FolnerDensity} implies \ref{item:UBD}, assume \ref{item:FolnerDensity}, and apply Lemma \ref{lem:MeanFromSequence} to find an invariant mean $m$ on $\ell^\infty(\Gamma)$ such that $m(1_A)=d_{\mb \Phi}(A)$.
		
		In the case where $\Gamma$ is countable, there are only countably many finite subsets of $\Gamma$, so the F{\o}lner net in \ref{item:FolnerDensity} may be taken to be a sequence.
	\end{proof}
	
	\begin{lemma}\label{lem:ThickEquivalents}   Let $A\subseteq \Gamma$.  The following are equivalent.
		
		\begin{enumerate}[label=(T.\arabic*)] 
			\item\label{item:isThick} $A$ is thick,
			\item $\Gamma\prec A$,
			\item for all finite $F\subseteq \Gamma$, $\bigcap_{\gamma \in F} (A-\gamma)\neq \varnothing$,
			\item $d^*(A)=1$,
			\item\label{item:MeanEquals1} $m(A)=1$ for some invariant mean $m$ on $\ell^\infty(\Gamma)$,
			\item $A\cap S\neq \varnothing$ for all syndetic sets $S\subseteq \Gamma$.
		\end{enumerate}
	\end{lemma}
	\begin{proof}
		(T.1)$\implies$(T.2) and (T.2)$\implies$(T.3) are both straightforward.  
		
		To prove (T.3)$\implies$(T.4), assume (T.3) and let $(\Phi_j)_{j\in I}$ be a F{\o}lner net for $\Gamma$.  For all $j\in I$, $\bigcap_{\gamma\in \Phi_j} (A-\gamma)\neq \varnothing$, so there is a $t_j\in \Gamma$ such that $t_j+\Phi_j\subseteq \Gamma$.  Then $\mathbf{\Phi}':=(t_j+\Phi_j)_{j\in I}$ is a F{\o}lner net such that  $d_{\mb \Phi}(A)=1$.  By Lemma \ref{lem:MeanFromSequence} there is an invariant mean $m$ such that $m(1_{A})=1$, so $d^{*}(A)=1$.
		
		Using Definition \ref{def:UBD}, (T.4)$\implies$(T.5) follows from the weak$^{*}$-compactness of $\mathcal M_{\tau}(\Gamma)$ and weak$^{*}$ continuity of the map $m\mapsto m(1_{A})$.
		
		To prove (T.5)$\implies$(T.6) let $m$ be an invariant mean with $m(1_A)=1$. If $S$ is syndetic then $\bigcup_{\gamma\in F} S+\gamma=\Gamma$ for some finite set $F$, so $m(S+F)=1$. Additivity of $m$ then implies $m(S)\geq 1/|F|$.  Since $m(A)=1$, we have $m(A\cap S)\geq 1/|F|$.  In particular $A\cap S$ is nonempty.
		
		To prove (T.6)$\implies$(T.1), let $A\subseteq \Gamma$ be such that $A\cap S\neq \varnothing$ for all syndetic sets $S$.  Suppose, to get a contradiction, that $A$ is not thick.  Then there is a finite set $F\subseteq \Gamma$ such that no translate of $F$ is contained in $A$.  In other words, $(\gamma +F)\cap A^{c}\neq \varnothing$ for all $\gamma \in \Gamma$.  It follows that for all $\gamma\in \Gamma$, $\gamma \in A^{c}-F$.  This means $A^{c}-F=\Gamma$, so $A^{c}$ is syndetic.  Now $A\cap A^{c}=\varnothing$, contradicting our assumption that $A\cap S\neq \varnothing$ for every syndetic set $S$.\end{proof}

	\begin{lemma}\label{lem:syndeticEquivalents}    Let $S\subseteq \Gamma$.  The following are equivalent.
		\begin{enumerate}
			\item[(S.1)] $S$ is syndetic.
			\item[(S.2)]  $S\cap T\neq \varnothing$ for all thick sets $T\subseteq \Gamma$.
			\item[(S.3)] $\Gamma \setminus S$ is not thick.
			\item[(S.4)] $d^*(\Gamma\setminus S)<1$.
		\end{enumerate}
	\end{lemma}
	\begin{proof}
		The implication (S.1)$\implies$(S.2) is the statement that every thick set has nonempty intersection with every syndetic set; this was proved in Lemma \ref{lem:ThickEquivalents}.
		
 (S.2)$\implies$(S.3) is obvious.
		
		To prove (S.3) implies (S.4), assume $\Gamma\setminus S$ is not thick.  Then Lemma \ref{lem:ThickEquivalents} implies $d^*(\Gamma\setminus S)<1$, so $1-d^*(\Gamma \setminus S)>0$.
		
		To prove (S.4) implies (S.1), we prove the contrapositive.  So assume that $S$ is not syndetic, meaning that for all finite $F\subseteq \Gamma$, we have $\bigcup_{\gamma\in F} S-\gamma \neq \Gamma$.  Writing $S^{c}$ for $\Gamma\setminus S$ and applying de Morgan's law, we get that for all finite $F\subseteq \Gamma$, $\bigcap_{\gamma\in F} S^c-\gamma \neq \varnothing$, meaning $S^c$ is thick.  Lemma \ref{lem:ThickEquivalents} then implies $d^*(S^c)=1$.	\end{proof}

\begin{lemma}\label{lem:PWSequivalents}
	Let $C\subseteq \Gamma$.  The following are equivalent.
	
	\begin{enumerate}
		\item[(PWS.1)]  There is a finite set $F\subseteq \Gamma$ such that $F+C$ is thick.
		\item[(PWS.2)]  There is a thick set $T\subseteq \Gamma$ and a syndetic set $S\subseteq \Gamma$ such that $S\cap T\subseteq C$.
		\item[(PWS.3)]  There is an invariant mean $m$ and a syndetic set $S\subseteq \Gamma$ such that $S\subset_{m} C$.
	\end{enumerate}
\end{lemma}

\begin{proof}
\noindent (PWS.1)$\implies$(PWS.2)  Assume $F$ is finite and $F+C$ is thick.   Let $T=C\cup (F+C)$, so that $T$ is also thick.  Let $S:= (C\cap T) \cup (\Gamma \setminus T)$.  

We claim that $S$ is syndetic.  Let $F'=F\cup \{0\}$. We will prove that $F'+S=\Gamma$.  Note that $C\subseteq S$, so $T=C\cup (F+C)\subseteq F'+S'$.  Since $0\in F'$, we have $0+(\Gamma\setminus T)=\Gamma\setminus T\subseteq F'+S$, as well.  This means $\Gamma = T\cup (\Gamma \setminus T)\subseteq F'+S$.  Thus $S$ is syndetic. Clearly $S\cap T \subseteq C$, so (PWS.2) is satisfied.

\noindent  (PWS.2)$\implies$(PWS.3)  Assuming $T$ is thick, $S$ is syndetic, and $S\cap T\subseteq C$, let $m$ be an invariant mean with $m(T)=1$.  Then $S\setminus C\subseteq S \setminus T$, so $m(S\setminus C)\leq m(S\setminus T)\leq m(\Gamma \setminus T)=0$.  This means $S\subset_{m} C$. 

\noindent (PWS.3)$\implies$(PWS.1) Assuming $S$ is syndetic and $S\subset_{m} C$, let $F$ be a finite set such that $F+S=\Gamma$.  Then $F+S\subset_{m} F+C$, so $\Gamma \subset_{m} F+C$, meaning $m(F+C)=1$.  By Lemma \ref{lem:ThickEquivalents} \ref{item:MeanEquals1}, this means $F+C$ is thick.  \end{proof}

	\subsection{Almost periodic functions and Bohr neighborhoods}\label{sec:BohrNhoods}

	Let $\mathbb T:=\mathbb R/\mathbb Z$ with the usual topology.  For $x\in \mathbb T$, let $\tilde{x}\in [0,1)$ be a coset representative of $x$ (i.e. $\tilde{x}+\mathbb Z = x$).  Let $\|x\|:=\min\{|\tilde{x}-n|:n\in \mathbb Z\}$.  Then $d_{\mathbb T}(x,y):=\|x-y\|$ is a translation invariant metric on $\mathbb T$.  For $\mb x = (x_1,\dots,x_d)\in \mathbb T^d$, let $\|\mb x\|:=\max_{j\leq d} \|x_j\|$.

	\begin{definition}
		Let $d\in \mathbb N$, let $\rho:\Gamma \to \mathbb T^d$ be a homomorphism, and let $\varepsilon>0$.  The \emph{Bohr$_0$-set} determined by these parameters is defined to be
		\[
		\Bohr(\rho;\varepsilon):=\{\gamma\in \Gamma: \|\rho(\gamma)\|<\varepsilon\}.
		\]
		In other words, $\Bohr(\rho;\varepsilon)$ is the preimage $\rho^{-1}(U)$, where $U:=\{\mb x\in \mathbb T^d: \|\mb x\|<\varepsilon\}$.  The \emph{rank} and \emph{radius} of $\Bohr(\rho;\varepsilon)$ are $d$ and $\varepsilon$, respectively.
		
		For a given $\gamma_0\in \Gamma$,  homomorphism $\rho:\Gamma\to \mathbb T^d$, and $\varepsilon>0$, the corresponding \emph{basic Bohr neighborhood of $\gamma_0$} is $\gamma_0 + \Bohr(\rho;\varepsilon)$; this is $\rho^{-1}(U)$, where $U=\{\mb x\in \mathbb T^d:\|\mb x-\rho(\gamma_{0})\|<\varepsilon\}$. We say that $U$ is a Bohr neighborhood of $\gamma_{0}$ if $U$ contains $\gamma_{0}+\Bohr(\rho;\varepsilon)$ for some  homomorphism $\rho:\Gamma\to \mathbb T^{d}$ and some $\varepsilon>0$.
	\end{definition}
	Since $\mathbb T^d$ is covered by at most $(\lfloor 1/\varepsilon\rfloor+1)^{d}$ translates of the open cube $\{\mb x:\|\mb x\|<\varepsilon\}$, we see that $\Gamma$ is covered by finitely many translates of $\Bohr(\rho;\varepsilon)$.  Thus $m(\Bohr(\rho;\varepsilon))\geq  (\lfloor 1/\varepsilon\rfloor + 1)^{-d}$ for every invariant mean $m$. Consequently every Bohr neighborhood is syndetic.

	\begin{lemma}\label{lem:BohrEquivalents}
		Let $A\subseteq \Gamma$ and $\gamma_{0}\in \Gamma$.  The following conditions are equivalent.
		
		\begin{enumerate}
			\item[(B.1)] $A$ contains a Bohr neighborhood of $\gamma_{0}$.
			
			\item[(B.2)] There is a compact abelian group $K$, a neighborhood $U$ of $0_K$, and a homomorphism $\rho:\Gamma\to K$ such that $\gamma_{0}+\rho^{-1}(U)\subseteq A$.
			
			\item[(B.3)] There is a uniformly almost periodic function $\phi:\Gamma\to [0,\infty)$ with $\phi(\gamma_{0})>0$ and $c<\phi(\gamma_{0})$ such that $\{\phi >c\}\subseteq A$. 
			
			\item[(B.4)] There is a neighborhood $U\subseteq b\Gamma$ of $0_{b\Gamma}$ such that $(\gamma_{0}+U)\cap \Gamma\subseteq A$.
		\end{enumerate}
	\end{lemma}
	
	\begin{proof}  We prove these equivalences in the special case where $\gamma_{0}=0$.  The general case is easy to obtain by translation.
		
		That (B.1)$\implies$(B.2) follows from the fact that $\mathbb T^d$ is a compact abelian group.
		
		To prove (B.2)$\implies$(B.3), let $K$ be a compact abelian group, let $\rho:\Gamma\to K$ be a homomorphism, and let $U$ be a neighborhood of $0_K$ in $K$.  Let $\tilde{f}:K\to [0,1]$ be a continuous function with $\tilde{f}(0_K)=1$ and $\tilde{f}(x)=0$ for all $x\notin U$.  Since $\tilde{f}$ is continuous, there is a trigonometric polynomial $q:=\sum_{j=1}^d c_j\psi_j$, where $\psi_j\in \widehat{K}$, such that $|\tilde{f}(x)-q(x)|<1/2$ for all $x\in K$.  Since $\tilde{f}$ is real-valued we can make $q$ real-valued (replace $q$ with $(q+\bar{q})/2$ if necessary).  Then $\chi_j:=\psi_j\circ \rho\in \widehat{\Gamma}$ for each $j$, so $\phi:=q\circ \rho$ is real-valued and uniformly almost periodic.  We also have $\{\gamma: \phi(\gamma)>1/2\}\subseteq \rho^{-1}(U)$, since $\phi(\gamma)>1/2$ implies $q(\rho(\gamma))>1/2$, which implies $\tilde{f}(\rho(\gamma))>0$.  Thus (B.3) is satisfied with $c=1/2$.
		
		To prove (B.3)$\implies$(B.1) let $\phi$ be a real-valued uniformly almost periodic function on $\Gamma$, let $c<\phi(0)$ and let $\psi:=\sum_{j=1}^d c_j\chi_j$ be a trigonometric polynomial such that $\|\phi - \psi\|_\infty<\delta$, where $\delta< \frac{1}{2}(\phi(0) - c)$.  Write $\chi_j(\gamma)$ as $\exp(2\pi i\rho_j(\gamma))$, where $\rho_j:\Gamma \to \mathbb T$ is a homomorphism, and let $\rho:\Gamma\to \mathbb T^d$ be $(\rho_1,\dots,\rho_d)$.  Choose $\varepsilon>0$ so that $|\exp(2\pi i t)-1|<\delta(1+\sum_{j=1}^d|c_j|)^{-1}$ whenever $|t|<\varepsilon$.

		Now if $\gamma\in \Bohr(\rho;\varepsilon)$ we have \[|\psi(\gamma)-\psi(0)|\leq \sum_{j=1}^d |c_j|(\exp(2\pi i\rho_j(\gamma))-1)\leq \delta \bigl(\sum_{j=1}^d |c_j|\bigr)\bigl(1+\sum_{j=1}^{d}|c_j|\bigr)^{-1}<\delta.\]
		Thus $|\phi(0)-\phi(\gamma)|\leq |\phi(0)-\psi(0)|+|\psi(0)-\phi(\gamma)|<2\delta$.  Our choice of $\delta$ then implies $\phi(\gamma)>c$.
		
		(B.3)$\implies$(B.4).  Let $\phi:\Gamma\to \mathbb R$ be uniformly almost periodic and $c<\phi(0)$ satisfy $\{\phi>c\}\subseteq A$. By Theorem \ref{th:BohrCompactificationAndAP}  there is a continuous $\tilde{\phi}:b\Gamma\to \mathbb R$ with $\tilde{\phi}|_\Gamma = \phi$.  Then $U:=\{\tilde{\phi}>c\}$ is a neighborhood of $0$ in $b\Gamma$ such that $U\cap \Gamma\subseteq A$.
		
		(B.4)$\implies$(B.2).  Note that (B.4) is the special case of (B.2) where $K=b\Gamma$ and $\rho$ is the embedding map.
	\end{proof}
	
		We say that $S\subseteq \Gamma$ is \emph{piecewise Bohr} if one of the (mutually equivalent) conditions in Lemma \ref{lem:PWBohrEquivalents} holds. 	
		  Using condition \ref{item:PWBdef}, we see that since every Bohr neighborhood is syndetic, every piecewise Bohr set is piecewise syndetic.   
	
	\begin{lemma}\label{lem:PWBohrEquivalents}
		Let $A\subseteq \Gamma$.  The following conditions are equivalent.
		
		\begin{enumerate}
			\item\label{item:PWBdef} there is a Bohr neighborhood $B$ of some $\gamma_{0}\in \Gamma$ and a thick set $T$ such that $B\cap T \subseteq S$.
			\item\label{item:PWBmean} There is an $m\in \mathcal M_{\tau}(\Gamma)$, a uniformly almost periodic function $\phi:\Gamma\to [0,\infty)$, and $c<\sup \phi$ such that $\{\phi>c\}\subset_m A$.
			\item\label{item:PWBprec} There is a Bohr$_0$-set $B$ such that $B\prec A$.
		\end{enumerate}
	\end{lemma}
	
	\begin{proof}
		\ref{item:PWBdef}$\implies$\ref{item:PWBmean}.  Assuming \ref{item:PWBdef}, there is a Bohr$_0$-set $B$, a $t\in \Gamma$, and a thick set $T$ such that $(t+B)\cap T\subseteq A$.  By Lemma \ref{lem:BohrEquivalents}, there is a uniformly almost periodic function $\phi_0:\Gamma\to \mathbb R$ and $c<\sup \phi_0$ such that $\{\phi_0>c\}\subseteq B$.  Since $T$ is thick, Lemma \ref{lem:ThickEquivalents} provides an invariant mean $m$ such that $m(T)=1$.  With this $m$ we have $\{\phi>c\}\subset_m A$, where $\phi$ is the translate of $\phi_0$ by $t$.

		\ref{item:PWBmean}$\implies$\ref{item:PWBprec}.  Let $\phi$ be almost periodic and $c<\sup \phi$ be such that $\{\phi>c\}\subset_m A$. Let $F\subseteq \{\phi>c\}$ be finite.  Let $c'>c$ so that $f(\gamma)>c'$ for all $\gamma\in F$.  By Theorem \ref{th:BohrCompactificationAndAP}, $\phi=\tilde{\phi}|_{\Gamma}$, where $\tilde{\phi}:b\Gamma\to \mathbb C$ is continuous.  Thus there is a neighborhood $U\subseteq b\Gamma$ of $0_{b\Gamma}$ such that $U+F\subseteq \{\tilde{\phi}>c\}$.  Letting $B=U\cap \Gamma$, we get that $B$ is a Bohr neighborhood of $0$ with $B+F\subseteq \{\phi>c\}$.  We then have $B\subseteq \bigcap_{\gamma\in F} \{\phi>c\}-\gamma$.  Since $\{\phi>c\}\subset_m A$, this implies $B\subset_m \bigcap_{\gamma \in F} A-\gamma$.  In particular, $\bigcap_{\gamma \in F} A-\gamma$ is nonempty, so $A$ contains a translate of $F$.  Since this holds for every finite $F\subseteq B$, we have $B\prec A$.

		\ref{item:PWBprec}$\implies$\ref{item:PWBdef}.  Fix a homomorphism $\rho:\Gamma\to \mathbb T^d$ and $\varepsilon>0$, and consider the corresponding Bohr$_0$-set $B=\Bohr(\rho;\varepsilon)$. Let $A$ be such that $B\prec A$.  For each finite $F\subseteq \Gamma$, choose $t_F\in \Gamma$ so that $(F\cap B)+t_F\subseteq A$.  Now $B+t_F=\{\gamma\in \Gamma: \|\rho(\gamma)-\rho(t_F)\|<\varepsilon\}$.  Let $I:=\{F\subseteq \Gamma: F \text{ is finite}\}$ be the directed set of finite subsets of $\Gamma$ with the usual preorder: $F<F'$ if $F\subseteq F'$.  Since $\mathbb T^d$ is compact, there is an $\alpha\in \mathbb T^d$ and a cofinal collection $\mathcal F\subseteq I$ such that $\|\rho(t_F)-\alpha\|<\varepsilon/4$ for every $F\in\mathcal F$.  Let $T=\bigcup_{F\in \mathcal F}F+t_F$;  $T$ is thick since $\mathcal F$ is cofinal in the partial order of finite subsets of $\Gamma$.
		Let $\gamma_0$ be such that $\|\alpha-\rho(\gamma_0)\|<\varepsilon/4$.  Let $B':=\{\gamma\in \Gamma:\|\rho(\gamma)-\rho(\gamma_0)\|<\varepsilon/2\}$.  We claim that $T\cap B'\subseteq A$.  To see this, let $\gamma\in T\cap B'$, so that $\|\rho(\gamma)-\rho(\gamma_0)\|<\varepsilon/2$ and for some $F\in \mathcal F$, $\gamma\in F+t_F$. 
		To prove that $\gamma\in A$, it suffices to show that $\gamma\in (F\cap B)+t_F$, meaning $\gamma-t_F\in F\cap B$.  We have $\gamma-t_F\in F$ already, we just have to show $\gamma-t_F\in B$.  To do so, write
		\[
		\|\rho(\gamma-t_F)\|=\|\rho(\gamma)-\rho(t_F)\|\leq \|\rho(\gamma)-\rho(\gamma_0)\|+\|\rho(\gamma_0)-\alpha\|+\|\alpha- \rho(t_F)\|<\varepsilon,
		\]
		so $\gamma-t_F\in B$.  \end{proof}

\section{Extreme points and maximizers}

\begin{lemma}\label{lem:UBDrealized}
	Let $K$ be a compact convex subset of a locally convex Hausdorff topological vector space $X$ and let $\lambda \in X^{*}$ be a continuous linear functional with $\lambda(K)=[a,b]\subseteq \mathbb R$.  Then there is an extreme point $k\in K$ such that $\lambda(k)=b$.
\end{lemma}

\begin{proof}
Continuity of $\lambda$  implies $K_{b}:=\lambda^{-1}(\{b\})\cap K$ is compact, and  linearity of $\lambda$ implies $K_{b}$ is convex.  By the Krein-Milman theorem, $K_{b}$ has extreme points.  Let $k\in K_{b}$ be such a point.  We will show that $k$ is an extreme point of $K$.  Write $k = cx+(1-c)y$, where $x$,  $y\in K$, $c\in [0,1]$.  We will show that $x=y=k$.  

Now $b=\lambda(k)=c\lambda(x)+(1-c)\lambda(y)$.  This implies $\lambda(x)=b$ and $\lambda(y)=b$, since $\lambda(K)\subseteq [a,b]$.  Thus $x$, $y\in K_{b}$, and we conclude that $x=y=k$, since $k$ is an extreme point of $K_{0}$.
\end{proof}

	\printbibliography


\end{document}
\providecommand{\MRhref}[2]{%
	\href{http://www.ams.org/mathscinet-getitem?mr=#1}{#2}
}
\providecommand{\href}[2]{#2}

\end{document}